\documentclass[draft]{conm-p-l}
\usepackage{amssymb}
\usepackage{graphicx}
\usepackage[mathscr]{eucal}
\copyrightinfo{2005}{American Mathematical Society}

\newtheorem{theorem}{Theorem}[section]
\newtheorem{lemma}[theorem]{Lemma}
\newtheorem{corollary}[theorem]{Corollary}
\newtheorem{hypothesis}[theorem]{Hypothesis}

\theoremstyle{definition}

\theoremstyle{remark}
\newtheorem{remark}[theorem]{Remark}
\numberwithin{equation}{section}

\newcommand{\bbR}{{\mathbb{R}}}
\newcommand{\bbP}{{\mathbb{P}}}

\newcommand{\bbC}{{\mathbb{C}}}

\newcommand{\calB}{{\mathcal B}}
\newcommand{\calC}{{\mathcal C}}

\newcommand{\calF}{{\mathcal F}}

\newcommand{\calH}{{\mathcal H}}

\newcommand{\calJ}{{\mathcal J}}
\newcommand{\calK}{{\mathcal K}}

\newcommand{\calU}{{\mathcal U}}

\newcommand{\scrD}{\mathscr{D}}
\newcommand{\scrH}{\mathscr{H}}

\newcommand{\no}{\notag}
\newcommand{\lb}{\label}
\newcommand{\f}{\frac}

\newcommand{\ol}{\overline}

\newcommand{\wti}{\widetilde}

\newcommand{\Oh}{O}
\newcommand{\oh}{o}

\newcommand{\loc}{\text{\rm{loc}}}

\newcommand{\dom}{\text{\rm{dom}}}

\newcommand{\supp}{\text{\rm{supp}}}
\newcommand{\AC}{\text{\rm{AC}}}
\newcommand{\bi}{\bibitem}
\newcommand{\hatt}{\widehat}
\newcommand{\beq}{\begin{equation}}
\newcommand{\eeq}{\end{equation}}
\newcommand{\ba}{\begin{align}}
\newcommand{\ea}{\end{align}}

\newcommand{\leftbrac}{\left(\negthickspace }
\newcommand{\rightbrac}{\negthickspace\right)}

\newcommand{\ga}{\gamma}
\newcommand{\la}{\lambda}

\newcommand{\ve}{\varepsilon}

\newcommand{\GopDhat}{G^{\opDhat}}
\newcommand{\GopDwti}{G^{\opDwti}}
\newcommand{\GopHhat}{G^{\opHhat}}
\newcommand{\GopHwti}{G^{\opHwti}}
\newcommand{\opDhat}{\hatt D}
\newcommand{\opDwti}{\wti D}
\newcommand{\opHhat}{\hatt H}
\newcommand{\opHwti}{\wti H}

\newcommand{\Ddag}{\scrD^{\dag}}
\newcommand{\DdagPsi}{\Psi^{\Ddag}}
\newcommand{\Dhat}{\hatt \scrD}
\newcommand{\Dhatdag}{\Dhat^{\dag}}
\newcommand{\DhatF}{\calF^{\Dhat}}
\newcommand{\DwtiF}{\calF^{\Dwti}}
\newcommand{\Dhatm}{m^{\Dhat}}
\newcommand{\DhatN}{N^{\Dhat}}
\newcommand{\DhatPhi}{\Phi^{\Dhat}}
\newcommand{\DhatPsi}{\Psi^{\Dhat}}
\newcommand{\DhatTheta}{\Theta^{\Dhat}}
\newcommand{\DPsi}{\Psi^{\scrD}}
\newcommand{\Dpsi}{\psi^{\scrD}}
\newcommand{\Dwti}{\wti \scrD}
\newcommand{\Dwtim}{m^{\Dwti}}
\newcommand{\DwtiN}{N^{\Dwti}}
\newcommand{\DwtiPhi}{\Phi^{\Dwti}}
\newcommand{\DwtiPsi}{\Psi^{\Dwti}}
\newcommand{\DwtiqPsi}{\Psi^{\Dwti_{q_0}}}
\newcommand{\DwtiTheta}{\Theta^{\Dwti}}

\newcommand{\Hhat}{\hatt \scrH}
\newcommand{\HhatF}{\calF^{\Hhat}}
\newcommand{\HwtiF}{\calF^{\Hwti}}
\newcommand{\Hhatm}{m^{\Hhat}}
\newcommand{\HhatN}{N^{\Hhat}}
\newcommand{\HhatPhi}{\Phi^{\Hhat}}

\newcommand{\HhatPsi}{\Psi^{\Hhat}}
\newcommand{\HPsi}{\Psi^{\scrH}}

\newcommand{\HhatT}{T^{\Hhat}}
\newcommand{\HwtiT}{T^{\Hwti}}
\newcommand{\HhatTheta}{\Theta^{\Hhat}}

\newcommand{\Hwti}{\wti \scrH}
\newcommand{\Hwtim}{m^{\Hwti}}
\newcommand{\HwtiN}{N^{\Hwti}}
\newcommand{\HwtiPhi}{\Phi^{\Hwti}}
\newcommand{\HwtiPsi}{\Psi^{\Hwti}}
\newcommand{\HwtiqPsi}{\Psi^{\Hwti_{q_0}}}
\newcommand{\HwtiTheta}{\Theta^{\Hwti}}

\newcommand{\mhat}{\hatt m}
\newcommand{\mwti}{\wti m}
\newcommand{\Mhat}{\hatt M}
\newcommand{\Mwti}{\wti M}
\newcommand{\Gammahat}{\hatt\Gamma}
\newcommand{\Gammawti}{\wti\Gamma}

\newcommand{\Rq}{\wti S_{q_{0}}}

\renewcommand{\Re}{\text{\rm Re}}
\renewcommand{\Im}{\text{\rm Im}}

\renewcommand{\ge}{\geqslant}
\renewcommand{\le}{\leqslant}

\newcounter{thmlistno}
\newenvironment{thmlist}
{\begin{list}
 {$(\roman{thmlistno})$}
 {\setlength{\itemindent}{0pt}
  \setlength{\labelwidth}{20pt}
  \setlength{\leftmargin}{23pt}
  \setlength{\itemsep}{0pt}
  \setlength{\topsep}{0pt}
  \setlength{\parskip}{0pt}
  \setlength{\parsep}{0pt}
  \setlength{\listparindent}{0pt}
  \usecounter{thmlistno} } }
{\end{list}}

\begin{document}

\title[On Self-adjoint and $J$-self-adjoint Dirac-type Operators]{On
Self-adjoint and $J$-self-adjoint Dirac-type \\ Operators: A Case
Study}

\author{Steve Clark}
\address{Department of Mathematics \& Statistics, University
of Missouri-Rolla, Rolla, MO 65409, USA}
\email{sclark@umr.edu}
\urladdr{http://www.umr.edu/\~{ }clark}
\thanks{The first author was supported in part by NSF Grant \#
DMS-0405526.}

\author{Fritz Gesztesy}
\address{Department of Mathematics,
University of Missouri, Columbia, MO 65211, USA}
\email{fritz@math.missouri.edu}
\urladdr{http://www.math.missouri.edu/personnel/faculty/gesztesyf.html}
\thanks{The second author was supported in part by NSF Grant \#DMS-0405528.}
\thanks{To appear in {\it Contemp. Math.}}




\begin{abstract}
We provide a comparative treatment of some aspects of spectral
theory for self-adjoint and non-self-adjoint (but
$J$-self-adjoint) Dirac-type operators connected with the
defocusing and focusing nonlinear Schr\"odinger equation, of relevance to
nonlinear optics.

In addition to a study of Dirac and Hamiltonian systems, we
also introduce the concept of Weyl--Titchmarsh half-line $m$-coefficients
(and $2\times 2$ matrix-valued $M$-matrices) in the non-self-adjoint
context and derive some of their basic properties. We conclude with an
illustrative  example showing that crossing spectral arcs in the
non-self-adjoint context imply the blowup of the norm of spectral projections
in the limit where the crossing point is approached.
\end{abstract}

\maketitle

\section{Introduction}

The principal part of this paper is devoted to a comparative study of
Dirac-type operators of the formally self-adjoint type
\begin{equation}
\Dhat= i\leftbrac
\begin{array}{cc}
             d/dx & -q(x) \\
             \ol{q(x)} & -d/dx
\end{array}\rightbrac , \quad x\in\bbR,    \lb{1.1}
\end{equation}
and the formally non-self-adjoint (but formally $J$-self-adjoint cf.\
\eqref{2.14a}) Dirac-type operators of the form
\begin{equation}
\Dwti = i\leftbrac
\begin{array}{cc}
             d/dx & -q(x) \\
             -\ol{q(x)} & -d/dx
\end{array}\rightbrac , \quad x\in\bbR,        \lb{1.2}
\end{equation}
where $q$ is locally integrable on $\bbR$.
Interest in these two particular Dirac-type operators stems from the fact
that both are intimately connected with applications to nonlinear optics. In
fact,  the differential expression $\Dhat$ gives rise to the Lax operator
of the {\it defocusing} nonlinear Schr\"odinger equation (${\rm NLS}_+$), while the
differential expression $\Dwti$ defines the Lax operator for the {\it
focusing} nonlinear Schr\"odinger equation (${\rm NLS}_-$). In appropriate
units, the propagation equation for a pulse envelope $q(x,t)$ in a monomode
optical fiber in the plane-wave limit neglecting loss is given by the
nonlinear Schr\"odinger equation
\begin{equation}
\text {NLS}_\pm (q):\,=iq_x \mp \tfrac{1}{2}q_{tt}+|q|^2q=0
\lb{1.3}
\end{equation}
(assuming weak nonlinearity of the medium and weak dispersion). The
focusing nonlinear Schr\"odinger equation admits a one-soliton
solution that propagates without change of shape and more generally
admits ``bright'' soliton solutions. The defocusing Schr\"odinger
equation shows a very different behavior since pulses undergo
enhanced broadening (to be used as optical pulse compression),
thereby yielding ``dark'' solitons. For pertinent general references
of this fascinating area we refer the reader, for instance, to
\cite{ADK93}, \cite{AA97}, \cite{DEGM88}, \cite{FT87}, \cite{GH03},
\cite{HK95}, \cite{HM03}, \cite{NM92}, \cite{Sh04}.

While typical applications to quantum mechanical problems in
connection with Schr\"odinger and Dirac equations require the study
of self-adjoint boundary value problems, many applications of
completely integrable systems most naturally lead to
non-self-adjoint Lax operators underlying the integrable system. The
prime example in this connection is the nonlinear Schr\"odinger
equation \eqref{1.3}. With this background in mind, we embarked upon
a more systematic study of the spectral properties of operator
realizations of \eqref{1.1} and especially, \eqref{1.2}, in
$L^2(\bbR)^2$.

There exists a large body of results on spectral and inverse spectral
theory of self-adjoint and non-self-adjoint Dirac-type operators,
especially, in the periodic and certain quasi-periodic cases (we refer,
e.g., to \cite[Ch.\ 5]{AK90}, \cite{CGHL04}--\cite{DM05a}, \cite{FL86},
\cite[Ch.\ 3]{GH03}, \cite{GT06}, \cite{GG93}--\cite{GKM98}, \cite{It81},
\cite{Kl05}, \cite{KS02}, \cite{Ko01}, \cite{Le91}, \cite{MA81},
\cite{Ma76}--\cite{Mi04}, \cite{Pr85}, \cite{Tk00},
\cite{Tk01}). It is impossible to refer to all relevant papers on the
subject, but a large list of references can be found in \cite{CG02}. In
this paper, however, we offer a different treatment focusing on a
comparative study of self-adjoint and non-self-adjoint (but
$J$-self-adjoint) Dirac-type operators with emphasis on
Weyl--Titchmarsh-type results. For basic results on $J$-self-adjoint 
operators we refer, for instance, to \cite[Sect.\ III.5]{EE89},
\cite[Sects.\ 21--24]{Gl65}, \cite{Kn81}, \cite{Ra85}, \cite{Zh59}. 
The Weyl--Titchmarsh $m$-coefficient was
first introduced for a class of $J$-self-adjoint Dirac-type operators with
bounded coefficients (and for the complex spectral parameter restricted to
a half-plane) in \cite{Sa90}. Additional results and further references
can be found in \cite{GKS98} and \cite{Sa05}. For a
general Weyl--Titchmarsh--Sims theory for singular non-self-adjoint
Hamiltonian systems we refer to \cite{BEP03}. Additional spectral results
and further references in the singular non-self-adjoint Hamiltonian system
case can be found in \cite{BM03}. 

In Section \ref{s2}, we begin by considering the general Dirac-type
expression
\begin{equation}
\scrD= \begin{pmatrix} i & 0 \\ 0 & -i \end{pmatrix}\f{d}{dx}+Q(x),
\quad Q=\begin{pmatrix} Q_{1,1} & Q_{1,2} \\ Q_{2,1} &
Q_{2,2} \end{pmatrix}\in L_{\loc}^1(\bbR)^{2\times 2},\; x\in\bbR.
\end{equation}
Introducing the conjugate linear operator acting upon $\bbC^2$, described by
\begin{equation}\label{calJ}
\calJ = \begin{pmatrix}0&1\\1&0 \end{pmatrix} \calC, \quad \calJ^2=I_2,
\end{equation}
where $\calC$ denotes the
operator of conjugation acting on $\bbC^2$ by
\begin{equation}\label{Conj}
\calC (a, b)^\top = (\ol a, \ol b)^\top, \quad a,b\in\bbC,
\end{equation}
with $(a, b)^\top$ denoting transposition of the vector $(a, b)$, we
show that $\scrD$ is formally $J$-self-adjoint (implying the same
property for $\Dhat$ and $\Dwti$). In particular, we show that
\begin{equation}
\calJ \scrD \calJ = \scrD^*,
\end{equation}
if and only if $Q_{1,1}=Q_{2,2}$ a.e.\ on $\bbR$.

Since Dirac-type operators are often studied in Hamiltonian form, we
also introduce the unitarily equivalent Hamiltonian form $\scrH$ of
$\scrD$ given by
\begin{equation}
\scrH=U\scrD U^{-1}= \begin{pmatrix} 0 & -1 \\ 1& 0 \end{pmatrix}\f{d}{dx} +B(x), \quad
B(x)=UQ(x)U^{-1}, \; x\in\bbR
\end{equation}
with the constant $2\times 2$ matrix $U$ given by \eqref{U}. Green's matrices are described
both for $\scrD$ and its unitarily equivalent Hamiltonian form $\scrH$.

Section \ref{s3} is devoted to a study of the maximally defined
$L^2(\bbR)^2$-realization $\hatt D$ of the special case $\hatt
\scrD$ in \eqref{1.1} and its unitarily equivalent Hamiltonian
version $\hatt H=U\hatt D U^{-1}$. $\hatt D$ (and hence $\hatt H$)
is known to be self-adjoint for all $q\in L^1_{\loc}(\bbR)$, (cf.
\cite{CG02}). We determine the Green's matrices of $\hatt H$ and
$\hatt D$ and recall some elements of the Weyl--Titchmarsh theory
associated with $\hatt H$. Due to the unitary equivalence of $\hatt
H$ and $\hatt D$, we show that the Weyl--Titchmarsh formalism for
$\hatt D$ can be set up in such a manner that the half-line
Weyl--Titchmarsh $m$-coefficients (and hence the $2\times 2$
matrix-valued full-line Weyl--Titchmarsh $M$-matrices) for $\hatt H$
and $\hatt D$ coincide. The latter appears to be new as
Weyl--Titchmarsh theory, to the best of our knowledge, is typically
formulated in connection with the Hamiltonian version $\hatt H$.
Moreover, we provide a streamlined derivation of the $2\times 2$
matrix-valued spectral functions of $\hatt H$ and $\hatt D$ starting
from the corresponding families of spectral projections. This
section is concluded with the simple constant coefficient example
$q(x)=q_0\in\bbC$ a.e.

Our final Section \ref{s4} then deals with a study of the maximally
defined $L^2(\bbR)^2$-realization $\wti D$ of the special case $\wti
\scrD$ in \eqref{1.2} and its unitarily equivalent Hamiltonian
version $\wti H=U\wti D U^{-1}$. $\wti D$ (and hence $\wti H$) is
known to be $J$-self-adjoint for all $q\in L^1_{\loc}(\bbR)$, (cf.
\cite{CG04}). We determine the Green's matrices of $\wti H$ and
$\wti D$, and develop some basic cornerstones of the analog of the
Weyl--Titchmarsh theory in the self-adjoint context of Section
\ref{s3} for the non-self-adjoint (but $J$-self-adjoint) operator
$\wti H$. Again, due to the unitary equivalence of $\wti H$ and
$\wti D$, we show that the Weyl--Titchmarsh formalism for $\wti D$
can be set up in such a manner that the half-line Weyl--Titchmarsh
$m$-coefficients (and hence the $2\times 2$ matrix-valued full-line
Weyl--Titchmarsh $M$-matrices) for $\wti H$ and $\wti D$ coincide.
In addition, we indicate the link between the spectral projections
of $\wti H$ (and hence of $\wti D$) and a $2\times 2$ matrix-valued
spectral function of $\wti H$ (and $\wti D$) determined from the
corresponding full-line Weyl--Titchmarsh $M$-matrix away from
spectral singularities of $\wti H$. This section also supplies the
illustrative constant coefficient example
$q(x)=q_0\in\bbC\backslash\{0\}$ a.e. In this case, the spectrum of
$\wti H$ consists of the real axis and the line segment from $-i|q_0|$ to
$+i|q_0|$ along the imaginary axis. In other words, this is
presumably the simplest differential operator with crossing spectral
arcs. We conclude this section with a proof of the fact that the
norm of the spectral projection in this example associated with an
interval of the type $(\lambda_1,\lambda_2)$,
$0<\lambda_1<\lambda_2$, blows up in the limit $\lambda_1\downarrow
0$, that is, when $\lambda_1$ approaches the crossing point
$\lambda=0$ of the spectral arcs of $\wti H$. The material developed
in this section represents the principal new results in this paper.

\section{A comparison of Dirac and Hamiltonian Systems}  \lb{s2}

\subsection{Dirac differential expressions}
Throughout this paper for a matrix $A$ with complex-valued entries,
$A^\top$ denotes the transposition of $A$; $\ol A$ denotes the
matrix with complex conjugate entries; and $A^*$ denotes the adjoint
matrix, that is, the conjugate transpose of $A$, $A^*=\ol A^\top$.
We will have occasion in our discussion to consider the following
$2\times 2$ matrices:
\begin{equation}\label{sigma}
\sigma_1=\begin{pmatrix}0&1\\1&0\end{pmatrix} \quad
I_2=\begin{pmatrix}1&0\\0&1\end{pmatrix}, \quad
\sigma_3=\begin{pmatrix}1&0\\0&-1\end{pmatrix}\quad
\sigma_4=\begin{pmatrix}0&1\\-1&0\end{pmatrix}.
\end{equation}
Moreover, we subsequently denote by $\sigma(A)$ and $\rho(A)$ the spectrum
and resolvent set of a closed densely defined linear operator $A$ in a
separable complex Hilbert space $\calH$.

We now consider \emph{whole-line Dirac differential expressions} of
the form
\begin{equation} \label{DE}
\scrD=i\sigma_3\dfrac{d}{dx} + Q(x), \quad  Q=\leftbrac
\begin{array}{cc}
             Q_{1,1} & Q_{1,2} \\
             Q_{2,1} & Q_{2,2}
\end{array}
\rightbrac\in L^1_{\loc}(\bbR)^{2\times 2},
\end{equation}
that is, $Q$ is a $2\times 2$ matrix with complex-valued entries that are
locally integrable on $\bbR$. In particular, we shall
be concerned with the formally self-adjoint differential expression that
arises when
\begin{equation}\label{DESA}
\Dhat=i\sigma_3\dfrac{d}{dx} + \hatt Q(x),\quad \hatt
Q=i\leftbrac\begin{array}{cc}
                      0 & -q \\
                      \ol{q} & 0
                    \end{array}\rightbrac\in L^1_{\loc}(\bbR)^{2\times 2},
\end{equation}
and the formally non-self-adjoint differential expression arising
when
\begin{equation}\label{DENSA}
\Dwti=i\sigma_3\dfrac{d}{dx} + \wti Q(x),\quad \wti
Q=i\leftbrac\begin{array}{cc}
                      0 & -q \\
                      -\ol{q} & 0
                    \end{array}\rightbrac\in L^1_{\loc}(\bbR)^{2\times 2}.
\end{equation}

By the  \emph{formal adjoint} of the differential expression
$\scrD$ given in \eqref{DE}, we shall mean the differential
expression $\scrD^*$, for which
\begin{equation}
\int_a^bdx\,\Psi(x)^*(\scrD \Phi)(x) =
\Psi(x)^*i\sigma_3\Phi(x)\big|_a^b +
\int_a^bdx\,\left(\scrD^*\Psi\right)(x)^* \Phi(x),
\end{equation}
for all $a, b \in \bbR$ and all
\begin{equation}
\Psi(x)=\begin{pmatrix}\psi_1(x)\\\psi_2(x) \end{pmatrix},\quad
\Phi(x)=\begin{pmatrix}\phi_1(x)\\\phi_2(x) \end{pmatrix},\quad
\psi_j,\ \phi_j \in\AC([a,b]), \; j=1,2,
\end{equation}
with $\AC([a,b])$ the set of absolutely continuous functions on $[a,b]$.
Hence, $\scrD^*$ is given by
\begin{equation}\label{2.4}
\scrD^*=i\sigma_3\dfrac{d}{dx} + Q(x)^*, \quad x\in\bbR.
\end{equation}
In particular, we note that $\Dhat^*=\Dhat$ while $\Dwti^*\ne
\Dwti$. Moreover, by the \emph{formal real adjoint} of the
differential expression $\scrD$ given by \eqref{DE}, we shall mean
the differential expression, $\Ddag$, where for all $a,b \in\bbR$,
\begin{equation}
\int_a^bdx\,\Psi (x)^\top(\scrD \Phi)(x) = \Psi (x)^\top
i\sigma_3\Phi(x)\big|_a^b + \int_a^bdx \left(\Ddag\Psi\right)
(x)^\top\Phi(x).
\end{equation}
Hence, $\Ddag$ is given by
\begin{equation}\label{2.6}
\Ddag=-i\sigma_3\dfrac{d}{dx} + Q(x)^\top.
\end{equation}

Associated with the Dirac differential expression \eqref{DE} is the
\emph{homogeneous Dirac system} given by
\begin{equation} \label{DS}
(\scrD\Psi)(z,x)=i\sigma_3\Psi'(z,x)+Q(x)\Psi(z,x)=z\Psi(z,x),
\end{equation}
for a.e.\ $x\in\bbR$, where $z$ plays the role of the spectral
parameter and
\begin{equation}
\Psi(z,x)=\begin{pmatrix}\psi_{1}(z,x)\\\psi_{2}(z,x)\end{pmatrix},\quad
\psi_j(z,\cdot)\in\AC_{\loc}(\bbR), \; j=1,2
\end{equation}
with $AC_{\loc}(\bbR)$ denoting the set of locally absolutely
continuous functions on $\bbR$.  By analogy, one obtains Dirac
systems associated with the differential expressions in
\eqref{DESA}, \eqref{DENSA}, \eqref{2.4}, and \eqref{2.6}. Solutions
of \eqref{DS} are said to be \emph{z-wave functions of} $\scrD$.

The Wronskian of two elements $F=(f_1 \; f_2)^\top, G=(g_1 \; g_2)^\top \in
C(\bbR)^2$ is defined as usual by
\begin{align}
\begin{split}
W(F(x),G(x))&=F(x)^\top\sigma_4 G(x)=f_1(x)g_2(x) - f_2(x) g_1(x) \\
&=\det\left(\begin{pmatrix} f_1(x) & g_1(x) \\ f_2(x) & g_2(x)
\end{pmatrix}\right), \quad x\in\bbR.
\end{split}
\end{align}

The differential expressions \eqref{DESA} and \eqref{DENSA}, which
will be the focus of our study, each exhibit the property of
\emph{formal $J$-self-adjointness}; a property that is manifest in
the following relations:
\begin{equation}
\calJ\Dhat\calJ=\Dhat^*=\Dhat,\quad \calJ\Dwti\calJ=\Dwti^*
\lb{2.14a}
\end{equation}
where $\calJ$ is defined in \eqref{calJ}, and where the equalities
hold a.e.\ on $\bbR$. While not all Dirac differential expressions
described in \eqref{DE} are formally $J$-self-adjoint, those which
are can be characterized as follows:

\begin{theorem}\label{tJSA}
Let $\scrD$ be the Dirac differential expression \eqref{DE}. Then
the following statements are equivalent:
\begin{thmlist}
\item $\scrD$ is formally $J$-self-adjoint:
$\calJ\scrD\calJ=\scrD^*$, where equality holds a.e.\ on $\bbR$.
\item  $Q_{1,1}=Q_{2,2}$ a.e.\ on $\bbR$ in the matrix $Q$ of the differential expression $\scrD$.
\item The Wronskian is a nonzero constant for any pair of linearly
independent $z$-wave functions of the Dirac system \eqref{DS}.
\end{thmlist}
\end{theorem}
\begin{proof}
The equivalence of statements $(i)$ and $(ii)$ follows from \eqref{2.4} and
the fact that
\begin{equation}
\calJ\scrD\calJ=i\sigma_3\frac{d}{dx} +
\begin{pmatrix}\ol Q_{2,2} & \ol Q_{2,1}\\ \ol Q_{1,2} & \ol
Q_{1,1}\end{pmatrix}.
\end{equation}
The equivalence of statements $(ii)$ and $(iii)$ follows from the
observation that if $\Psi_j(z,x)$, $j=1,2$ represent two independent
$z$-wave functions of the Dirac system \eqref{DS}, then
\begin{align}
\frac{d}{dx}W(\Psi_1(z,x),\Psi_2(z,x))
&=i[Q_{1,1}(x)-Q_{2,2}(x)]W(\Psi_1(z,x),\Psi_2(z,x)),
\end{align}
and hence
\begin{equation}\lb{Wronskian}
W(\Psi_1(z,x),\Psi_2(z,x))=W(\Psi_1(z,0),\Psi_2(z,0))
\exp{\bigg(i\int_0^x ds\, [Q_{1,1}(s)-Q_{2,2}(s)]\bigg)}.
\end{equation}
\end{proof}
           From \eqref{Wronskian}, we obtain
\begin{corollary}
The Wronskian has nonzero constant magnitude for any pair of
independent $z$-wave functions of the Dirac system \eqref{DS} if
and only if
\begin{equation}
            \Im [Q_{1,1}(x)-Q_{2,2}(x)] = 0 \, \text{ for a.e.\ $x\in\bbR$.}
\end{equation}
\end{corollary}
In light of Theorem~\ref{tJSA} and the earlier observation that our
study will focus upon the two examples of $J$-self-adjoint
differential expressions provided by $\Dhat$ and $\Dwti$, we make the
following hypothesis for the remainder of this paper:
\begin{hypothesis}\label{hJSA}
We assume that the Dirac differential expression $\scrD$ given in
\eqref{DE} is formally $J$-self-adjoint, that is,
$Q_{1,1}(x)=Q_{2,2}(x)$ holds for a.e.\ $x\in\bbR$.
\end{hypothesis}

\subsection{Green's matrices and Dirac operators}\label{sDMDO}

Assuming the existence of a whole-line Green's matrix for a
$J$-self-adjoint Dirac system \eqref{DS}, we can associate a
Dirac operator $D$ on $\bbR$ in the following way: Let $f\in
L^2(\bbR)^2$, assume $\rho\subset\bbC$ is open and nonempty, and consider
the inhomogeneous Dirac system given by
\begin{equation}\label{IHDS}
(\scrD\Psi)(z,x)=i\sigma_3\Psi'(z,x)+Q(x)\Psi(z,x)=z\Psi(z,x) +
f(x), \quad z\in\rho.
\end{equation}
If $G^{D}(z,x,x')$, $z\in\rho$, $x,x'\in\bbR$, denotes the unique Green's
matrix associated with \eqref{DS}, then \eqref{IHDS} has a unique solution,
$\Psi(z,\cdot)\in L^2(\bbR)^2\cap \AC_{\loc}(\bbR)^2$ given by
\begin{equation}
\Psi(z,x) = \int_{-\infty}^{\infty}\, dx' \, G^{D}(z,x,x')f(x'), \quad
z\in\rho, \; x\in\bbR.
\end{equation}
The Dirac operator $D$ in $L^2(\bbR)^2$ associated with \eqref{DS}
is then defined by
\begin{equation}
((D-z)^{-1}f)(x)= \int_{-\infty}^{\infty}\, dx' \, G^{D}(z,x,x')
f(x'), \quad z\in\rho, \; f\in L^2(\bbR)^2.
\end{equation}
In terms of the differential expression \eqref{DE}, $D$ is
explicitly defined by
\begin{align}
\begin{split}
D&=i\sigma_3 \dfrac{d}{dx} +Q,  \\
\dom(D)&=\big\{\Psi\in L^2(\bbR)^2\cap \AC_{\loc}(\bbR)^2 \,\big|\,
\scrD\Psi\in L^2(\bbR)^2 \big\}. \label{DO}
\end{split}
\end{align}

\begin{remark}
Construction of a unique whole-line Green's matrix for \eqref{DS} in
association with the operator $D$ is equivalent to the existence of unique
(up to constant multiples) Weyl--Titchmarsh-type solutions
$\DPsi_{\pm}(z,\cdot)\in L^2([0,\pm\infty))^2$, $z\in\rho$ of \eqref{DS}.
Such solutions are known to exist for Dirac systems associated with
\eqref{DESA} (cf.\ \cite{CG02}), and \eqref{DENSA} (cf.\ \cite{CGHL04}).
Hence by the construction above, one can describe the operator $\opDhat$ in
association with \eqref{DESA} and the operator $\opDwti$ in association
with \eqref{DENSA}.  As special cases of \eqref{DE} satisfying
Hypothesis~\ref{hJSA}, both
$\Dhat$ and $\Dwti$ are formally $J$-self-adjoint differential
expressions. Moreover, it has been proved in \cite{CG02} that $\hatt D$,
maximally defined as in \eqref{DO}, is self-adjoint,
\begin{equation}
\hatt D=\hatt D^*.
\end{equation}
In addition, it was shown in \cite{CGHL04} that $\wti D$, maximally
defined as in \eqref{DO}, is $J$-self-adjoint,
\begin{equation}
\calJ \wti D\calJ=\wti D^*.
\end{equation}
In the self-adjoint context \eqref{DESA}, the existence of unique
Weyl--Titchmarsh-type solutions is of course equivalent to the limit point
case of $\Dhat$ at $\pm\infty$. In the context of \eqref{DENSA}, the
existence of unique Weyl--Titchmarsh-type solutions, or equivalently, the
existence of a unique Green's function, is then the proper  analog of the
limit point case in this non-self-adjoint situation.
\end{remark}

In the next Lemma, and under the presumption of
the existence of half-line square integrable solutions, we describe
the whole-line Green's matrix for the Dirac system \eqref{DS} in
association with the operator $D$ defined in \eqref{DO}.

\begin{lemma}\label{l2.1}
Let $\rho \subset\bbC$ be open and nonempty. Suppose that for all
$z\in\rho$,
$ \DPsi_{\pm}(z,\cdot)=
\leftbrac\begin{array}{c}
             \Dpsi_{\pm,1}(z,\cdot) \\
             \Dpsi_{\pm,2}(z,\cdot)
\end{array}\rightbrac
\in L^2([0,\pm\infty))^{2}$ represent a basis of solutions of the
Dirac system given by \eqref{DS}. Then the whole-line Green's matrix
for this system is given by
\begin{align}\label{2.9}
G^{D}(z,x,x')&=C(z) \begin{cases} \DPsi_-(z,x) \DPsi_+(z,x')^\top
\sigma_1, & x<x', \\[15pt]
\DPsi_+(z,x) \DPsi_-(z,x')^\top\sigma_1,& x>x'
\end{cases}\\ \notag
&=C(z) \begin{cases}
\leftbrac\begin{array}{cc}
             \Dpsi_{-,1}(z,x)\Dpsi_{+,2}(z,x') &
\Dpsi_{-,1}(z,x)\Dpsi_{+,1}(z,x') \\
             \Dpsi_{-,2}(z,x)\Dpsi_{+,2}(z,x') &
\Dpsi_{-,2}(z,x)\Dpsi_{+,1}(z,x')
\end{array}\rightbrac,& x<x', \\[15pt]
\leftbrac\begin{array}{cc}
             \Dpsi_{+,1}(z,x)\Dpsi_{-,2}(z,x') &
\Dpsi_{+,1}(z,x)\Dpsi_{-,1}(z,x') \\
             \Dpsi_{+,2}(z,x)\Dpsi_{-,2}(z,x') &
\Dpsi_{+,2}(z,x)\Dpsi_{-,1}(z,x')
\end{array}\rightbrac,& x>x',
\end{cases} \\
& \hspace*{9.8cm} z\in\rho, \no
\end{align}
where
\begin{align}
\begin{split}
C(z)&=-i[W(\DPsi_+(z,x),\DPsi_-(z,x))]^{-1} \\
&=-i\big[\Dpsi_{+,1}(z,x)
\Dpsi_{-,2}(z,x)-\Dpsi_{-,1}(z,x)\Dpsi_{+,2}(z,x)\big]^{-1}
\end{split}
\end{align}
is constant with respect to $x\in\bbR$.
\end{lemma}
\begin{proof}
Note the following unitary equivalence of differential expressions
associated with \eqref{DE} and \eqref{2.6}:
\begin{equation}\label{2.14}
\scrD - z = \sigma_1 (\Ddag - z )\sigma_1.
\end{equation}
As a consequence, if
\begin{equation}\label{2.10}
\DdagPsi_{\pm}(z,\cdot)= \sigma_1 \DPsi_{\pm}(z,\cdot)\in
L^2([0,\pm\infty))^{2}, \quad z\in\rho,
\end{equation}
then $\DdagPsi_{\pm}(z,x)$ represent half-line square integrable
solutions of the associated real adjoint system
\begin{equation}
(\Ddag\Psi)(z,x)=-i\sigma_3\Psi'(z,x)+Q(x)^\top \Psi(z,x) =
z\Psi(z,x).
\end{equation}
Hence, the Green's matrix ansatz given in \eqref{2.9} can be written
as
\begin{equation}\label{2.12}
G^{D}(z,x,x')= \begin{cases}
C\DPsi_-(z,x)\DdagPsi_+(z,x')^\top, & x<x',\\
C\DPsi_+(z,x)\DdagPsi_-(z,x')^\top, & x>x'.
\end{cases}
\end{equation}

To verify the ansatz,
let $f \in L^2(\bbR)^{2}$ and note that
\begin{align}
&(\scrD-z)\int_{-\infty}^\infty\, dx' \, G^{D}(z,x,x')f(x')\\
&=(\scrD-z)\int_{-\infty}^x\, dx' \,
G^{D}(z,x,x')f(x')+(\scrD-z)\int_x^\infty\,
dx' \, G^{D}(z,x,x')f(x')\notag\\
&=iC\sigma_3\big[\DPsi_+(z,x)\DdagPsi_-(z,x)^\top -
\DPsi_-(z,x)\DdagPsi_+(z,x)^\top \big]f(x) \notag
\end{align}
with the last equality following from the fact that
$(\scrD-z)G^D(z,x,x')=0$ for $x\ne x'$. Given \eqref{2.10}, we see
that
\begin{align}
\begin{split}
&(\scrD-z)\int_{-\infty}^\infty\, dx' \, G^{D}(z,x,x')f(x')\\
&=iC\leftbrac
\begin{array}{cc}
           ( \Dpsi_{+,1}\Dpsi_{-,2}-\Dpsi_{-,1}\Dpsi_{+,2})(z,x) & 0 \\
            0 & (\Dpsi_{-,2}\Dpsi_{+,1}-\Dpsi_{+,2}\Dpsi_{-,1})(z,x)
\end{array}\rightbrac f(x).
\end{split}
\end{align}
Given that the Wronksian
\begin{align}
\begin{split}
W(\DPsi_+(z,x),\DPsi_-(z,x))&=\DPsi_-(z,x)^\top J\DPsi_+(z,x) \\
&=\Dpsi_{+,1}(z,x)\Dpsi_{-,2}(z,x)-\Dpsi_{-,1}(z,x)\Dpsi_{+,2}(z,x)
\end{split}
\end{align}
is a nonzero constant for $x\in\bbR$, we obtain
\begin{equation}
f(x)=(\scrD-z)\int_{-\infty}^\infty\, dx' \, G^{D}(z,x,x')f(x')
\end{equation}
when $C=-i(\Dpsi_{+,1}\Dpsi_{-,2}-\Dpsi_{-,1}\Dpsi_{+,2})^{-1}$.
\end{proof}

\subsection{Hamiltonian Systems and Green's matrices}

Associated with the whole-line formally $J$-self-adjoint Dirac
differential expression \eqref{DE} is the unitarily equivalent
differential expression in Hamiltonian form given by
\begin{equation} \label{HE}
\scrH=U\scrD U^{-1}=-\sigma_4\dfrac{d}{dx} + B(x),  \quad B\in
L^1_{\loc}(\bbR)^{2\times 2},
\end{equation}
in terms of the unitary matrix
\begin{equation}\label{U}
U=\frac12\leftbrac
\begin{array}{rr}
             -1+i & -1+i \\
             1+i & -1-i
\end{array}
\rightbrac, \quad U^*=U^{-1}.
\end{equation}
In particular, we note that
\begin{equation}
B=UQU^{-1}= \dfrac12
\leftbrac\begin{array}{cc}
             Q_{2,1}+Q_{1,2}+2Q_{1,1} & i(Q_{2,1}-Q_{1,2}) \\
             i(Q_{2,1}-Q_{1,2}) & -Q_{2,1}-Q_{1,2}+2Q_{1,1}
\end{array}\rightbrac.
\end{equation}
With $U\in\bbC^{2\times 2}$ defined in \eqref{U}, we note also that
\begin{equation}
i\sigma_4=U\sigma_3 U^{-1}, \quad U\sigma_1 U^\top =-iI_2,
\end{equation}
and observe that the property of formal $J$-self-adjointness for
$\scrD$ is now manifest in the unitarily equivalent Hamiltonian
differential expression $\scrH$ by the following relationships:
\begin{equation}\label{HJSA}
(i\calC)\scrH(i\calC)=\scrH^* ,
\end{equation}
where $\calC$ again represents the conjugation operator \eqref{Conj} acting
on $\bbC^2$.

We note that in association with \eqref{DESA}, one obtains the
unitarily equivalent formally self-adjoint differential expression
$\Hhat$ given by
\begin{equation}\label{HESA}
\Hhat=U\Dhat U^{-1}=-\sigma_4\dfrac{d}{dx} + \hatt B(x),\quad \hatt
B=U\hatt QU^{-1}=\leftbrac\begin{array}{rr}
             \Im (q) & -\Re (q) \\
             -\Re (q) & -\Im (q)
\end{array}\rightbrac,
\end{equation}
while in association with \eqref{DENSA} one obtains the unitarily
equivalent formally non-self-adjoint differential expression
$\Hwti$ given by
\begin{equation}\label{HENSA}
\Hwti =U\Dwti U^{-1}= -\sigma_4\dfrac{d}{dx} + \wti B(x),\quad \wti
B=U\wti QU^{-1}=i\leftbrac\begin{array}{rr}
             -\Re(q) & -\Im(q) \\
             -\Im(q) & \Re(q)
\end{array}\rightbrac.
\end{equation}
As special cases of \eqref{HE}, both $\Hhat$ and $\Hwti$
satisfy the relationship given in \eqref{HJSA}.

In addition to providing unitarily equivalent differential
expressions, the unitary matrix $U$ exhibits another notable
feature: It preserves the Wronskian.

\begin{lemma}\label{l2.2}
Let $\eta=\begin{pmatrix}\eta_1\\ \eta_2\end{pmatrix},\,
\xi=\begin{pmatrix}\xi_1\\ \xi_2\end{pmatrix}\in\bbC^2$. Then
\begin{equation}
\det (\eta\ \, \xi)= \det(U \eta\ \, U\xi)=\det(U^{-1}\eta\ \, U^{-1}\xi).
\end{equation}
\end{lemma}
\begin{proof}
We note that $\det (\eta\ \, \xi)= \eta^\top \sigma_4\xi$, where
$\sigma_4\in\bbC^{2\times 2}$ is defined in \eqref{sigma}.
Similarly, we note that $\det(U \eta\ \, U\xi)= \eta^\top U^\top
\sigma_4U\xi$, and that $\det(U^{-1}\eta\ \,  U^{-1}\xi)=\eta^\top
(U^{-1})^\top JU^{-1}\xi$.  Upon verifying that $\sigma_4=U^\top
\sigma_4U=(U^{-1})^\top \sigma_4U^{-1}$, the result follows.
\end{proof}

Through the unitary equivalence in \eqref{HE} of the differential
expressions $\scrD$ and $\scrH$, we can define an operator $H$, in
association with the homogeneous Hamiltonian system given by
\begin{equation} \label{HS}
(\scrH\Psi)(z,x)=-\sigma_4\Psi'(z,x)+B(x)\Psi(z,x)=z\Psi(z,x),
\end{equation}
for a.e.\ $x\in\bbR$, where $z$ plays the role of the spectral
parameter, and where
\begin{equation}
\Psi(z,x)=\begin{pmatrix}\psi_{1}(z,x)\\\psi_{2}(z,x)\end{pmatrix},\quad
\psi_j(z,\cdot)\in \AC_{\loc}(\bbR), \; j=1,2.
\end{equation}
Namely,
\begin{align} \label{HO}
H&=UD U^{-1}= -\sigma_4 \dfrac{d}{dx} +B\\
\dom(H)&=\big\{\Psi\in L^2(\bbR)^2 \cap \AC_{\loc}(\bbR)^2 \big|\,
\scrH\Psi\in L^2(\bbR)^2 \}.
\end{align}

The presumptive existence of half-line square integrable
solutions of the Dirac system \eqref{DS} yields the existence of
half-line square integrable solutions of the associated Hamiltonian
system \eqref{HS} by
\begin{equation}
\HPsi_\pm(z,x)=U\DPsi_\pm(z,x).
\end{equation}
Then, by  Lemma~\ref{l2.1}, we  obtain a description of the
whole-line Green's matrix for the Hamiltonian system \eqref{HS}.

\begin{lemma}\label{l2.3}
Let $\rho \subset\bbC$ be open and nonempty. Suppose that for all
$z\in\rho$, $\DPsi_{\pm}(z,\cdot)
\in L^2([0,\pm\infty))^{2} $ represent a basis of solutions of the
Dirac system given by \eqref{DS}. Then, with $U$ given in \eqref{U},
$\HPsi_\pm(z,\cdot)
=U\DPsi_\pm(z,x)\in L^2([0,\pm\infty))^{2}$,
represent a basis of solutions of the Hamiltonian
system given by \eqref{HS}, and the whole-line Green's matrix for
the Hamiltonian system \eqref{HS} is given by
\begin{align}
G^{H}(z,x,x')&= UG^{D}(z,x,x')U^{-1} \no \\
&=K(z) \HPsi_\mp(z,x)\HPsi_\pm(z,x')^\top , \quad x\lessgtr
x', \; z\in\rho,
\end{align}
where
\begin{equation}
K(z)=[W(\HPsi_+(z,x),\HPsi_-(z,x))]^{-1}
=[W(\DPsi(z,x)_+,\DPsi_-(z,x))]^{-1}
\end{equation}
is constant with respect to $x\in\bbR$.
\end{lemma}
\begin{proof}
Let $z\in\rho$. By the unitary equivalence of $\scrD$ and $\scrH$
seen in \eqref{HO}, and with $U$ defined in \eqref{U} and
$\sigma_1\in\bbR^{2\times 2}$ defined in \eqref{sigma}, it follows
that
\begin{align}
G^{H}(z,x,x')&= UG^{D}(z,x,x')U^{-1}\\
&=C U\DPsi_\mp(z,x)\DPsi_\pm(z,x')^\top \sigma_1U^{-1}, \quad x\lessgtr
x', \no \\
&=C \HPsi_\mp(z,x)\HPsi_\pm(z,x')^\top (U^{-1})^\top \sigma_1 U^{-1},
\quad x\lessgtr x', \no \\
&=iC\HPsi_\mp(z,x)\HPsi_\pm(z,x')^\top , \quad x\lessgtr x', \no
\end{align}
where by Lemmas \ref{l2.1} and \ref{l2.2},
\begin{align}
\begin{split}
C&=-i[W(\DPsi_+,\DPsi_-)]^{-1}\\
&=-i[W(U^{-1}\HPsi_+,U^{-1}\HPsi_-)]^{-1} =-i[W(\HPsi_+,\HPsi_-)]^{-1}.
\end{split}
\end{align}
\end{proof}

\section{Self-adjoint Dirac and Hamiltonian Systems}  \lb{s3}

As developed in the previous section, the Dirac operator $D$
defined in \eqref{DO} corresponding to the Dirac system
\eqref{DS}, is unitarily equivalent to the operator $H$ in
\eqref{HO} associated with the Hamiltonian system \eqref{HS}. In
this section, we focus upon self-adjoint realizations for each of
these operators, specifically, the operator $\opDhat$, maximally defined by
\eqref{DO} associated with the special case of \eqref{DE}
given by \eqref{DESA}, and the operator $\opHhat$ maximally defined by
\eqref{HO} corresponding to the special case of \eqref{HE}
given by \eqref{HESA}.

\subsection{Weyl--Titchmarsh coefficients.}\label{sWTCSA}

Let $\HhatN(z,\pm\infty)$ and $\DhatN(z,\pm\infty)$, $z\in\bbC$, denote the
spaces  defined for the differential expressions $\Hhat$ and $\Dhat$,
respectively,  by
\begin{align} \lb{3.1}
\HhatN(z,\pm\infty)&= \big\{ \Psi\in L^2([0,\pm\infty))^2\,\big|\, (\Hhat
-z)\Psi =0 \big\},\\
\DhatN(z,\pm\infty)&= \big\{ \Psi\in L^2([0,\pm\infty))^2\,\big|\, (\Dhat
-z)\Psi =0 \big\}.\label{3.1b}
\end{align}
By \cite[Lemma 2.15]{CG02},
\begin{equation}\lb{3.2a}
\dim \big(\HhatN(z,\pm\infty)\big)=1, \quad
z\in\bbC\backslash\bbR,
\end{equation}
and hence by the unitary equivalence
given in \eqref{HE},
\begin{equation}\lb{3.2b}
\dim\big(\DhatN(z,\pm\infty)\big)=1,\quad z\in\bbC\backslash\bbR.
\end{equation}
In particular, one has the following result.

\begin{theorem} [\cite{CG02}, \cite{LM03}, \cite{We71}] \lb{t3.1.1}
The operator $\hatt H$, maximally defined in \eqref{HO}, is
self-adjoint,
\begin{equation}
\hatt H=\hatt H^*   \lb{3.4}
\end{equation}
and the operator $\hatt D$, maximally defined in \eqref{DO}, is
self-adjoint,
\begin{equation}
\hatt D =\hatt D^*.  \lb{3.5}
\end{equation}
Moreover, $\hatt H$ and $\hatt D$ are unitarily equivalent,
\begin{equation}
\hatt H=U\hatt D U^{-1}.  \lb{3.6a}
\end{equation}
\end{theorem}
\begin{proof}
Equation \eqref{3.4} has been proven in \cite{CG02}. The rest follows from
the unitary equivalence \eqref{HE} via the constant unitary matrix $U$.
\end{proof}
Self-adjoint half-line operators associated with
the differential expressions $\Hhat$ and $\Dhat$ are defined by
\begin{align}
\opHhat_\pm(\alpha)&=-\sigma_4\dfrac{d}{dx} + \hatt B, \no \\
\dom(\opHhat_\pm(\alpha))&= \big\{ \Psi\in
L^2([0,\pm\infty))^2\,\big|\, \Psi\in\AC_\loc([0,\pm\infty))^2,  \lb{3.8a}
\\ & \hspace*{2.05cm} \alpha\Psi(0)=0,\ \Hhat\Psi\in L^2([0,\pm\infty))^2
\big\}, \no
\end{align}
where $\alpha=(\cos(\theta),\sin(\theta))$, $\theta\in [0,2\pi)$, and by
\begin{align}
\opDhat_\pm(\beta)&=i\sigma_3\dfrac{d}{dx} + \hatt Q, \no \\
\dom(\opDhat_\pm(\beta))&= \big\{ \Psi\in
L^2([0,\pm\infty))^2\,\big|\, \Psi\in\AC_\loc([0,\pm\infty))^2,  \lb{3.9a}
\\ & \hspace*{2.1cm} \beta\Psi(0)=0,\ \Dhat\Psi\in L^2([0,\pm\infty))^2
\big\},
\no
\end{align}
where
\begin{equation}
\beta=\alpha U=[(-1+i)/2](e^{-i\theta},e^{i\theta}),
\quad \theta\in [0,2\pi).
\end{equation}

$\opHhat_\pm(\alpha)$ is unitarily equivalent
to $\opDhat_\pm(\beta)$, given \eqref{HESA} and the fact that the unitary
$2\times 2$ matrix $U$ naturally defines a unitary mapping of
$L^2([0,\pm\infty))^2$ onto itself, again for simplicity denoted by $U$,
which maps $\dom(\opDhat_\pm(\beta))$ onto $\dom(\opHhat_\pm(\alpha))$. The
later fact can be seen by noting that
\begin{equation}
0=\beta\DhatPsi(z,0)=\beta U^{-1}U\DhatPsi(z,0)=\alpha\HhatPsi(z,0).
\end{equation}
Thus,
\begin{equation}
\opHhat_\pm(\alpha)=U \opDhat_\pm(\beta) U^{-1}, \quad \beta=\alpha U.
\end{equation}
In passing, we note that \eqref{3.2a}
and \eqref{3.2b} prove that both $\opHhat_\pm(\alpha)$ and
$\opDhat_\pm(\beta)$ are in the limit point case at $\pm\infty$.

Next, let a fundamental system of solutions of the self-adjoint
Hamiltonian system $\Hhat \Psi=z\Psi$ be given by
\begin{equation} \label{FSHhat}
\HhatTheta(z,\cdot,\alpha),\
\HhatPhi(z,\cdot,\alpha)\in \AC_\loc(\bbR)^2,
\quad z\in\bbC
\end{equation}
such that
\begin{equation}
\HhatTheta(z,0,\alpha)=\alpha^*, \quad
\HhatPhi(z,0,\alpha)=-\sigma_4\alpha^*,
\end{equation}
where $\alpha\in\bbC^2$ and where
\begin{equation}
\alpha\alpha^*=1,\quad \alpha\sigma_4\alpha^*=0.  \lb{3.15a}
\end{equation}
By Theorem~\ref{tJSA}, $W(\HhatTheta(z,x,\alpha), \HhatPhi(z,x,\alpha))$ is constant for $x\in\bbR$.
If in addition, we require that
\begin{equation}\label{FSHhatc}
W(\HhatTheta(z,0,\alpha), \HhatPhi(z,0,\alpha)) = -\det(\alpha^*\
\sigma_4\alpha^*)=1,
\end{equation}
then it can be shown that $\alpha\in\bbR^2$. Thus, \eqref{3.8a} and
\eqref{3.9a} yield the only self-adjoint half-line operators consistent
with \eqref{3.15a} and \eqref{FSHhatc}. Hence, for the remainder of this
section, we let $\alpha=\alpha(\theta)=(\cos(\theta), \sin(\theta))$, for
$\theta\in [0,2\pi)$, and let
\begin{equation}\label{FSHhatd}
\HhatTheta(z,0,\alpha)=(\cos(\theta), \sin(\theta))^\top,\quad
\HhatPhi(z,0,\alpha)=(-\sin(\theta),\cos(\theta))^\top.
\end{equation}

In \cite{CG02}, it is shown that $\HhatPhi(z,\cdot,\alpha)\notin
L^2([0,\pm\infty))^2$ for $z\in\bbC\backslash\bbR$. Then, as a
consequence of \eqref{3.2a}, let $\Hhatm_\pm (z,\alpha)$ denote the
half-line Weyl--Titchmarsh coefficients; that is, the unique coefficients
such that
\begin{equation}\label{WSHhat}
\HhatPsi_\pm(z,\cdot,\alpha) =\HhatTheta(z,\cdot,\alpha) +
\Hhatm_\pm(z,\alpha)\HhatPhi(z,\cdot,\alpha)\in L^2([0,\pm\infty))^2, \quad
z\in\bbC\backslash\bbR.
\end{equation}

A corresponding development for the self-adjoint Dirac system
$\Dhat\Psi(z,x)=z\Psi(z,x)$ begins with its fundamental system of
solutions
\begin{equation} \label{FSDhat}
\DhatTheta(z,\cdot,\beta),\ \DhatPhi(z,\cdot,\beta)\in \AC_\loc(\bbR)^2,
\quad z\in\bbC
\end{equation}
for $\beta=\alpha U$, where $\alpha=(\cos(\theta),\sin(\theta))$,
$\theta\in [0,2\pi)$ (cf.\ \eqref{FSHhatd}),
$U\in\bbC^{2\times 2}$ is given in \eqref{U}, and hence,
\begin{equation}\label{beta}
\beta=\alpha U=[(-1+i)/2](e^{-i\theta},
e^{i\theta}), \quad\theta\in [0,2\pi).
\end{equation}
Specifically, for $\theta \in [0,2\pi) $, let
\begin{align}
\DhatTheta(z,0,\beta)&=i \sigma_1 \beta^\top
=-[(1+i)/2] (e^{i\theta},e^{-i\theta})^\top ,   \lb{3.21}\\
\DhatPhi(z,0,\beta)&=i\sigma_3\beta^*=[(1-i)/2]
(e^{i\theta},-e^{-i\theta})^\top.  \lb{3.22}
\end{align}
In particular, we see that
\begin{equation}\label{3.6}
\DhatTheta(z,x,\beta)=U^{-1}\HhatTheta(z,x,\alpha), \quad
\DhatPhi(z,x,\beta)=U^{-1}\HhatPhi(z,x,\alpha).
\end{equation}
As a consequence, $\DhatPhi(z,\cdot,\beta)\notin L^2([0,\pm\infty))^2$. As
before, given \eqref{3.1b}, let $\Dhatm_\pm(z,\beta)$ denote
the unique coefficients such that
\begin{equation}\label{WSDhat}
\DhatPsi_\pm(z,\cdot,\beta)
=\DhatTheta(z,\cdot,\beta) + \Dhatm_\pm(z,\beta)\DhatPhi(z,\cdot,\beta)\in
L^2([0,\pm\infty))^2, \quad z\in\bbC\backslash\bbR.
\end{equation}
In summary, we have the following result.

\begin{lemma}\label{l3.1}
Let $\alpha=(\cos(\theta),\sin(\theta))$, $\theta\in [0,2\pi)$, and let
$\beta = \alpha U$ with $U$ defined in \eqref{U}. Let $\HhatTheta$,
$\HhatPhi$ represent the fundamental system of solutions of the
Hamiltonian system $\Hhat\Psi=z\Psi$ satisfying \eqref{FSHhatd}, and
let $\DhatTheta$, $\DhatPhi$ represent the fundamental system of
solutions of the Dirac system $\Dhat\Psi=z\Psi$ satisfying
\eqref{3.21}, and \eqref{3.22}. Then, with $\HhatPsi_\pm$
defined in
\eqref{WSHhat} and with $\DhatPsi_\pm$ defined in \eqref{WSDhat}, one
infers that
\begin{equation}
\DhatPsi_\pm(z,x,\beta)=U^{-1}\HhatPsi_\pm(z,x,\alpha),
\quad z\in\bbC\backslash\bbR, \; x\in\bbR,
\end{equation}
and in particular, that
\begin{equation}\label{3.8}
\Dhatm_\pm(z,\beta)=\Hhatm_\pm(z,\alpha), \quad z\in\bbC\backslash\bbR.
\end{equation}
\end{lemma}
\begin{proof}
By \eqref{3.6}
\begin{align}
\begin{split}
U^{-1}\HhatPsi_\pm(z,\cdot,\alpha)&=U^{-1}\HhatTheta(z,\cdot,\alpha) +
\Hhatm_\pm(z,\alpha) U^{-1}\HhatPhi(z,\cdot,\alpha)\\
&=\DhatTheta(z,\cdot,\beta) +
\Hhatm_\pm(z,\alpha)
\DhatPhi(z,\cdot,\beta)\in L^2([0,\pm\infty))^2.
\end{split}
\end{align}
By the uniqueness of the representation for the combination given in
\eqref{WSDhat}, equation \eqref{3.8} follows.
\end{proof}

\begin{remark}\label{r3.2}
In light of Lemma~\ref{l3.1},  in the future we shall represent both
$\Hhatm_\pm(z,\alpha)$ and $\Dhatm_\pm(z,\beta)$ by $\mhat_{\pm}(z,\gamma)$,
where it is understood that $\beta=\alpha U$. Here
\begin{equation}
\gamma = \begin{cases} \text{represents $\alpha$ in the context of
$\hatt\scrH$}, \\
\text{represents $\beta=\alpha U$ in the context of
$\hatt\scrD$}, \end{cases}
\end{equation}
and we keep this convention in similar contexts in the
following.
\end{remark}

Of course, $\pm \hatt m_{\pm}(\cdot,\gamma)$ are well-known to be Herglotz
functions (i.e., analytic functions mapping the open complex upper
half-plane into itself) and
\begin{equation}
\hatt m_{\pm}(\cdot,\gamma) \, \text{ are analytic on } \,
\rho(\hatt H_{\pm}(\alpha))=\rho(\hatt D_{\pm}(\beta)).
\end{equation}

\subsection{Green's matrices}

Before describing Green's matrices for self-adjoint Hamiltonian and
Dirac systems, we introduce two matrices. First, for the fundamental
system of solutions defined in \eqref{FSHhat} and satisfying \eqref{FSHhatd},
let $\HhatF(z,\cdot,\alpha)$ denote the associated fundamental matrix given
by
\begin{equation}
\HhatF(z,x,\alpha)=\big(\HhatTheta(z,x,\alpha) \ \
\HhatPhi(z,x,\alpha)\big), \quad z\in\bbC, \; x\in\bbR.
\end{equation}
Next, we introduce the matrix $\Gammahat(z,\gamma)$ (we recall the meaning
of $\gamma$ as introduced in Remark \ref{r3.2}), where
\begin{equation}\label{Gammahat}
\Gammahat(z,\gamma)=\begin{pmatrix}
\frac1{\mhat_-(z,\gamma) - \mhat_+(z,\gamma)}&
\frac{\mhat_-(z,\gamma)}{\mhat_-(z,\gamma) - \mhat_+(z,\gamma)}\\[6pt]
\frac{ \mhat_+(z,\gamma)}{\mhat_-(z,\gamma) - \mhat_+(z,\gamma)}&
\frac{\mhat_-(z,\gamma)\mhat_+(z,\gamma)}{\mhat_-(z,\gamma)
- \mhat_+(z,\gamma)}
\end{pmatrix}, \quad z\in\bbC\backslash\bbR.
\end{equation}
Then, as a consequence of Lemma~\ref{l2.3} we obtain the following result.

\begin{lemma}\label{l3.3}
With $\HhatPsi_\pm(z,\cdot,\alpha)$ representing the half-line
Weyl--Titchmarsh solutions defined in \eqref{WSHhat} for the Hamiltonian
system $\Hhat\Psi=z\Psi$, the associated whole-line Green's
matrix associated with the operator $\opHhat$ is given by
\begin{align}
\GopHhat(z,x,x')&= K(z,\alpha)
\HhatPsi_\mp(z,x,\alpha)\HhatPsi_\pm(z,x',\alpha)^\top ,\quad
x\lessgtr x', \label{3.12a} \\
&=\begin{cases}
\HhatF(z,x,\alpha)\Gammahat(z,\gamma)^\top \HhatF(z,x'\alpha)^\top ,
& x<x',\\
\HhatF(z,x,\alpha)\Gammahat(z,\gamma)\HhatF(z,x'\alpha)^\top , & x>x',
\end{cases}  \quad z\in\bbC\backslash\bbR, \label{3.12b}
\end{align}
where
\begin{equation}
K(z,\alpha)=\big[W(\HhatPsi_+(z,x,\alpha),
\HhatPsi_-(z,x,\alpha))\big]^{-1}\big|_{x=0}=
[\mhat_-(z,\gamma)-\mhat_+(z,\gamma)]^{-1}.
\end{equation}
\end{lemma}
\begin{proof}
Equation \eqref{3.12a} follows from Lemma~\ref{l2.3} for the
operator $\opHhat$ defined by \eqref{HO}, but in association with the
special case of \eqref{HE} given by \eqref{HESA}. Moreover, it
follows that $K=[W(\HhatPsi_+,\HhatPsi_-)]^{-1}$. Then, by
\eqref{WSHhat} one notes that
\begin{equation}
W(\HhatPsi_+,\HhatPsi_-) =\begin{pmatrix} 1\\ \mhat_-
\end{pmatrix}^\top
\begin{pmatrix}
(\HhatTheta)^\top J\HhatTheta&(\HhatTheta)^\top J\HhatPhi\\[2pt]
(\HhatPhi)^\top J\HhatTheta& (\HhatPhi)^\top J\HhatPhi
\end{pmatrix}
\begin{pmatrix}
1\\\mhat_+
\end{pmatrix}.
\end{equation}
However, $\eta^\top J\eta=0$, and $\eta^\top J\xi=-\xi^\top J\eta$
for every $\eta,\ \xi\in\bbC^2$. As a consequence,
\begin{align}
W(\HhatPsi_+,\HhatPsi_-)&=[\mhat_- -
\mhat_+](\HhatPhi)^\top J\HhatTheta \no \\
&=[\mhat_- - \mhat_+]W(\HhatTheta, \HhatPhi) \no \\
&=\mhat_- - \mhat_+,
\end{align}
where the last equality follows from the normalization \eqref{FSHhatc}.

The description of $\GopHhat(z,x,x')$ given in \eqref{3.12b} follows
from \eqref{WSHhat}, \eqref{3.12a},  and the fact that
\begin{equation}\label{3.15}
\HhatPsi_\pm (z,\cdot,\alpha)= \HhatF(z,\cdot,\alpha)\begin{pmatrix}1\\
\mhat_\pm(z,\alpha)
\end{pmatrix}.
\end{equation}
\end{proof}

Following as an immediate consequence of the unitary equivalence of
$\opHhat$ and $\opDhat$, together with Lemmas~\ref{l2.1},
\ref{l2.3}, and \ref{l3.3}, one infers the following fact.

\begin{lemma}\label{l3.4}
With $\HhatPsi_\pm(z,\cdot,\alpha)$ and $\DhatPsi_\pm(z,\cdot,\beta)$
defined in \eqref{WSHhat} and \eqref{WSDhat}, the whole-line
Green's matrix for the self-adjoint Dirac system $\Dhat\Psi
=z\Psi$ is given by
\begin{align}
\GopDhat (z,x,x')&=C(z,\beta) \DhatPsi_\mp(z,x,\beta)\Psi_\pm^{\Dhatdag}
(z,x',\beta)^\top ,\quad x\lessgtr x', \no \\
&=C(z,\beta)\DhatPsi_\mp(z,x,\beta)\DhatPsi_\pm(z,x',\beta)^\top
\sigma_1,\quad
x\lessgtr x',\label{3.16b} \\
&= \begin{cases}
-i\DhatF(z,x,\beta)\Gammahat(z,\gamma)^\top \DhatF(z,x',\beta)^\top ,
& x<x',\\
-i\DhatF(z,x,\beta)\Gammahat(z,\gamma)\DhatF(z,x',\beta)^\top , & x>x',
\end{cases} \quad z\in\bbC\backslash\bbR, \no
\end{align}
where
\begin{align}
C(z,\beta)&=-i\big[W(\DhatPsi_+(z,x,\beta),
\DhatPsi_-(z,x,\beta))\big]^{-1}\big|_{x=0}
\no \\
&=-i\big[W(\HhatPsi_+(z,x,\alpha),
\HhatPsi_-(z,x,\alpha))\big]^{-1}\big|_{x=0} \no \\
&=-i[\mhat_-(z,\gamma)-\mhat_+(z,\gamma)]^{-1}.
\end{align}
Here $\DhatF(z,\cdot,\beta)$ is the fundamental matrix of solutions of the
Dirac system $\Dhat \Psi =z\Psi$ given by
\begin{equation}
\DhatF(z,x,\beta)=\big(\DhatTheta(z,x,\beta) \ \,
\DhatPhi(z,x,\beta)\big), \quad z\in\bbC, \; x\in\bbR,
\end{equation}
and $\DhatTheta(z,\cdot,\beta)$, $\DhatPhi(z,\cdot,\beta)$ represent the
fundamental system of solutions defined in \eqref{FSDhat} for the
self-adjoint Dirac system.
\end{lemma}

Of course, the Green's matrices \eqref{3.12a} and \eqref{3.16b} extend to
analytic $2\times 2$ matrix-valued functions with respect to $z\in
\rho(\hatt H) = \rho(\hatt D)$.

\subsection{Spectral matrices} \lb{GMSA}

In preparation for the description of the spectral matrix associated
with the operator $\opHhat$, we now introduce two matrices and a
transformation.

We denote by $\hatt M(z,\gamma)\in\bbC^{2\times 2}$,
$z\in\bbC\backslash\bbR$, the whole-line Weyl--Titchmarsh $M$-function
of the operator $\opHhat$ defined in
\eqref{HO} in association with the special case given by
\eqref{HESA}, namely,
\begin{align}\label{Mhat}
\Mhat(z,\gamma)&=\big(
\Mhat_{\ell,\ell'}(z,\gamma)\big)_{\ell,\ell'=0,1}
=\frac12[\hatt\Gamma(z,\gamma)+\hatt\Gamma(z,\gamma)^\top]
= \Gamma(z,\gamma) +\f{1}{2}\begin{pmatrix} 0 & -1 \\ 1 & 0 \end{pmatrix}
\no \\
&=\begin{pmatrix} \frac1{\mhat_-(z,\gamma) - \mhat_+(z,\gamma)}&
\frac12\frac{\mhat_-(z,\gamma) +
\mhat_+(z,\gamma)}{\mhat_-(z,\gamma) -
\mhat_+(z,\gamma)}\\[6pt]
\frac12\frac{\mhat_-(z,\gamma) + \mhat_+(z,\gamma)}{\mhat_-(z,\gamma) -
\mhat_+(z,\gamma)}&
\frac{\mhat_-(z,\gamma)\mhat_+(z,\gamma)}{\mhat_-(z,\gamma) -
\mhat_+(z,\gamma)}
\end{pmatrix}, \quad z\in\bbC\backslash\bbR,
\end{align}
where by Remark \ref{r3.2},
$\mhat_\pm(z,\gamma)=\Hhatm_\pm(z,\alpha)=\Dhatm_\pm(z,\beta)$ for
$\beta=\alpha U$, and $\alpha=(\cos(\theta),\sin(\theta))$, $\theta\in
[0,2\pi)$. Again, \eqref{Mhat} extends to an analytic $2\times 2$
matrix-valued function with respect to $z\in\rho(\hatt H) =
\rho(\hatt D)$.

From $\Im \big(\Mhat(z,\gamma)\big)
=(2i)^{-1}\big[\Mhat(z,\gamma) -
\Mhat(z,\gamma)^*\big]$ one infers that
\begin{equation}\label{3.20}
\Im \big(\Mhat(z,\gamma)\big)=\begin{pmatrix}
\Im\left(\frac1{\mhat_-(z,\gamma) - \mhat_+(z,\gamma)}\right)&
\frac12\Im\left(\frac{\mhat_-(z,\gamma) +
\mhat_+(z,\gamma)}{\mhat_-(z,\gamma) -
\mhat_+(z,\gamma)}\right)\\[6pt]
\frac12\Im\left(\frac{\mhat_-(z,\gamma) +
\mhat_+(z,\gamma)}{\mhat_-(z,\gamma) -
\mhat_+(z,\gamma)}\right)&
\Im\left(\frac{\mhat_-(z,\gamma)\mhat_+(z,\gamma)}{\mhat_-(z,\gamma)
-\mhat_+(z,\gamma)}\right)
\end{pmatrix}.
\end{equation}
Associated with $\hatt M (z,\gamma)$ we introduce the measure
$d\hatt\Omega(\lambda,\gamma)$ by
\begin{equation}
\hatt\Omega((\lambda_1.\lambda_2],\gamma)=\f{1}{\pi}
\lim_{\delta\downarrow 0}\lim_{\varepsilon\downarrow 0}
\int_{\lambda_1+\delta}^{\lambda_2+\delta} d\lambda \,
\Im\big(\hatt M(\lambda+i\varepsilon,\gamma)\big), \quad \lambda_j\in\bbR,
\; j=1,2, \; \lambda_1<\lambda_2,  \lb{3.42}
\end{equation}
and use the abbreviation
\begin{equation}
\big( \HhatT_0(\alpha) f\big)(\lambda) = \int_\bbR dx \,
\HhatF(\lambda,x,\alpha)^\top f(x), \quad \lambda \in \bbR, \; f\in
C^\infty_0(\bbR)^2.
\end{equation}
Henceforth we also abbreviate the scalar product in $L^2((a,b))^2$ by
$\langle \cdot, \cdot \rangle_{L^2((a,b))^2}$ (chosen to be linear in the
second place), where $-\infty\leq a<b\leq \infty$.

\begin{theorem}\label{tWSHhat}
Let $\big\{E_{\opHhat}(\lambda)\big\}_{\lambda\in\bbR}$ denote the spectral
family associated with the operator $\opHhat$. Then, for $f,\, g\in
C_0^\infty(\bbR)^2$ and $\lambda_1<\lambda_2$,
\begin{equation}\label{3.13}
\langle f, E_{\opHhat}((\lambda_1,\lambda_2])g\rangle_{L^2(\bbR)^2}=
\int_{(\lambda_1,\lambda_2]}\left(\big( \HhatT_0(\alpha)
f\big)(\lambda)\right)^*d\hatt\Omega(\lambda,\gamma) \big(
\HhatT_0(\alpha) g\big)(\lambda).
\end{equation}
\end{theorem}
\begin{proof}
For simplicity we will suppress the $\alpha$ (resp., $\gamma$) dependence
of all quantities involved in this proof. We follow the strategy of proof
employed in connection with one-dimensional Schr\"odinger operators in
\cite{GZ05} (see also \cite{HS98}).  Then, by Stone's formula (cf.\
\cite[p.\ 191]{We80}),
\begin{align}
&\langle f, E_{\opHhat}((\lambda_1,\lambda_2])g\rangle_{L^2(\bbR)^2} \no \\
& \quad =\lim_{\delta \downarrow 0}\lim_{\varepsilon \downarrow 0}
\frac1{2\pi i}\int_{\lambda_1+\delta}^{\lambda_2+\delta}d\lambda \,
\big\langle f,\
\big[(\opHhat-(\lambda+i\varepsilon))^{-1}-
(\opHhat-(\lambda-i\varepsilon))^{-1}\big]g \big\rangle_{L^2(\bbR)^2}
\notag \\ & \quad =\lim_{\delta \downarrow 0}\lim_{\varepsilon
\downarrow 0}\frac1{2\pi
i}\int_{\lambda_1+\delta}^{\lambda_2+\delta}d\lambda
\int_{\bbR}dx\int_{\bbR}dx'\big\{
f(x)^*\GopHhat(\lambda+i\varepsilon, x,x')g(x')  \no \\
& \hspace*{6cm} - f(x)^*\GopHhat(\lambda-i\varepsilon,x,x')g(x')\big\}.
\lb{3.46}
\end{align}
Using the fact that
$\mhat_\pm(\lambda-i\varepsilon)=\overline{\mhat_\pm(\lambda
+i\varepsilon)}$, one concludes that
$\Gammahat(\lambda-i\varepsilon)
=\overline{\Gammahat(\lambda+i\varepsilon)}$,
where $\Gammahat(z)$ is defined in \eqref{Gammahat}.
Consequently, using the description of $\GopHhat(z,x,x')$ given
in \eqref{3.12b}, \eqref{3.46} implies
\begin{align}
&\langle f, E_{\opHhat}((\lambda_1,\lambda_2])g\rangle_{L^2(\bbR)^2}
=\lim_{\delta \downarrow 0}\lim_{\varepsilon \downarrow 0}\frac1
{2\pi i}\int_{\lambda_1+\delta}^{\lambda_2+\delta}d\lambda\int_{\bbR}dx
\no \\
&\hspace{1.2cm} \times \bigg\{\int_{-\infty}^xdx'
f(x)^*\HhatF(\lambda,x)\Big[\Gammahat(\lambda+i\varepsilon)
- \overline{\Gammahat(\lambda+i\varepsilon)} \,
\Big]\HhatF(\lambda,x')^\top g(x') \no \\
& \hspace*{1.8cm} +\int_x^{\infty}dx'f(x)^*\HhatF(\lambda,x)
\Big[\Gammahat(\lambda+i\varepsilon)^\top
- \overline{\Gammahat(\lambda+i\varepsilon)^\top}
\, \Big]\HhatF(\lambda,x')^\top g(x')\bigg\}  \no  \\
& \hspace*{.9cm }
=\lim_{\delta \downarrow 0}\lim_{\varepsilon \downarrow 0}\frac1
{\pi}\int_{\lambda_1+\delta}^{\lambda_2+\delta}d\lambda\int_{\bbR}dx
\no \\
&\hspace{2cm} \times \bigg\{\int_{-\infty}^xdx'
f(x)^*\HhatF(\lambda,x) \Im\big(\hatt M(\lambda+i\varepsilon)\big)
\HhatF(\lambda,x')^\top g(x') \no \\
& \hspace*{2.7cm} +\int_x^{\infty}dx'f(x)^*\HhatF(\lambda,x)
\Im\big(\hatt M(\lambda+i\varepsilon)\big)
\HhatF(\lambda,x')^\top g(x')\bigg\},  \lb{3.48a}
\end{align}
since
\begin{equation}
\hatt\Gamma(\lambda+i\varepsilon)-\ol{\hatt\Gamma(\lambda+i\varepsilon)}=
\hatt\Gamma(\lambda+i\varepsilon)^\top-\ol{\hatt\Gamma
(\lambda+i\varepsilon)^\top}=2i \,
\Im\big(\hatt M(\lambda+i\varepsilon)\big),
\quad \lambda\in\bbR, \; \varepsilon>0.
\end{equation}
To arrive at equation \eqref{3.48a} we used the fact that for fixed
$x\in\bbR$, $\HhatF(z,x)$ is entire with respect to $z$, that
$\HhatF(\lambda,x)$ is real-valued for $\la\in\bbR$, that
$\HhatF(\lambda,x)^\top \in AC_{\loc}(\bbR)^2$, and hence that
\begin{equation}
\HhatF(\lambda\pm i\ve,x) \underset{\varepsilon\downarrow 0}{=}
\HhatF(\lambda,x) \pm i\varepsilon (d/dz)\HhatF(z,x)|_{z=\lambda} +
\Oh(\ve^2), \lb{2.85a}
\end{equation}
with $\Oh(\varepsilon^2)$ being uniform with respect to
$(\lambda,x)$ as long as $(\lambda,x)$ vary in compact subsets of
$\bbR^2$. Moreover, we used that
\begin{align}
\begin{split}
&\varepsilon|\hatt M_{\ell,\ell'}(\lambda+i\varepsilon,\ga)|\leq
C(\lambda_1,\lambda_2,\varepsilon_0), \quad \lambda\in
[\lambda_1,\lambda_2], \; 0<\varepsilon\leq\varepsilon_0, \;
\ell,\ell'=0,1,  \\
&\varepsilon |\Re(\hatt M_{\ell,\ell'}(\lambda+i\varepsilon,\ga))|
\underset{\varepsilon\downarrow 0}{=}\oh(1), \quad \lambda\in\bbR,
\; \ell,\ell'=0,1,
\end{split}
\end{align}
which follow from the properties of Herglotz functions since
$\hatt M_{\ell,\ell}$, $\ell=0,1$, are Herglotz and
$\hatt M_{0,1}=\hatt M_{1,0}$ have Herglotz-type representations by
decomposing the associated complex measure $d\hatt\Omega_{0,1}$ into
$d\hatt\Omega_{0,1}=d(\omega_1-\omega_2)+id(\omega_3-\omega_4)$, with
$d\omega_k$, $k=1,\dots,4$, nonnegative measures. Finally, we also used
(for $\lambda\in\bbR$, $\varepsilon>0$)
\begin{align}
\begin{split}
\hatt\Gamma(\lambda+i\varepsilon)+\ol{\hatt\Gamma(\lambda+i\varepsilon)}
&=2\, \Re\big(\hatt M(\lambda+i\varepsilon)\big)+\begin{pmatrix}
0 & 1 \\ -1 & 0 \end{pmatrix}, \\
\hatt\Gamma(\lambda+i\varepsilon)^\top
+\ol{\hatt\Gamma(\lambda+i\varepsilon)^\top}
&=2\, \Re\big(\hatt M(\lambda+i\varepsilon)\big)+\begin{pmatrix}
0 & -1 \\ 1 & 0 \end{pmatrix}. \\
\end{split}
\end{align}
Thus,
\begin{align}\label{3.31}
&\langle f, E_{\opHhat}((\lambda_1,\lambda_2])g\rangle_{L^2(\bbR)^2}  \\
&\quad = \int_{(\lambda_1,\lambda_2]} \int_{\bbR}dx
\int_\bbR dx' f(x)^*\HhatF(\lambda,x) d\hatt \Omega (\lambda)
\HhatF(\lambda,x')^\top g(x').\notag
\end{align}
Equation \eqref{3.13} then follows from the fact that $\HhatF(\bar
z,\cdot)=\overline{\HhatF(z,\cdot)}$, $z\in\bbC$, and hence that
\begin{equation}
\left(\big( \HhatT_0 f\big)(\lambda)\right)^*= \int_\bbR dx f(x)^*
\overline{\HhatF(\lambda,x)}= \int_\bbR dx f(x)^*\HhatF(\lambda,x),
\quad \lambda\in\bbR.
\end{equation}
\end{proof}

The proof of Theorem \ref{tWSHhat} shows that $\HhatT_0 (\alpha)$
represents a linear operator (denoted by the same symbol),
\begin{equation}
 \HhatT_0 (\alpha)\colon\begin{cases} C_0^\infty(\bbR)^2 \to L^2(\bbR;
d\hatt\Omega(\lambda,\gamma)) \\
\hspace*{1.06cm} f \mapsto  \HhatT_0(\alpha) f = \int_\bbR dx \,
\HhatF(\cdot,x,\alpha)^\top f(x).
\end{cases}
\end{equation}
(For some subtleties of $L^2$-spaces with matrix-valued measures we
refer to the discussion in \cite{GZ05} and the references cited
therein.) Moreover, as recently discussed in the analogous context
of Schr\"odinger operators in \cite{GZ05}, $ \HhatT_0(\alpha)$
extends to a bounded operator from $L^2(\bbR)^2$ to $L^2(\bbR;
d\hatt\Omega(\lambda,\gamma))$, which we denote by $\HhatT(\alpha)$.
This then immediately leads to the following extension of Theorem
\ref{tWSHhat}.

\begin{theorem} \lb{t3.7}
Let $\big\{E_{\opHhat}(\lambda)\big\}_{\lambda\in\bbR}$ denote the spectral
family associated with the operator $\opHhat$. Then, for $f,\, g\in
L^2 (\bbR)^2$ and $\lambda_1<\lambda_2$,
\begin{equation}\label{3.52}
\langle f, E_{\opHhat}((\lambda_1,\lambda_2])g\rangle_{L^2(\bbR)^2} =
\int_{(\lambda_1,\lambda_2]}\left(\big(\HhatT(\alpha)
f\big)(\lambda)\right)^*d\hatt\Omega(\lambda,\gamma)
\big(\HhatT(\alpha) g\big)(\lambda).
\end{equation}
\end{theorem}

As a corollary, we obtain the corresponding result for the operator
$\opDhat$.

\begin{corollary}
Let $\big\{E_{\opDhat}(\lambda)\big\}_{\lambda\in\bbR}$ denote the spectral
family associated with the operator $\opDhat$. Then, for $f,\, g\in
L^2(\bbR)^2$ and $\lambda_1<\lambda_2$,
\begin{equation}
\langle f, E_{\opDhat}((\lambda_1,\lambda_2])g\rangle_{L^2(\bbR)^2} =
\int_{(\lambda_1,\lambda_2]}\left(\big(\HhatT(\alpha)
Uf\big)(\lambda)\right)^*
d\hatt\Omega(\lambda,\gamma)\big(\HhatT(\alpha) Ug\big)(\lambda).
\end{equation}
\end{corollary}
\begin{proof}
This follows immediately from the observation that
\begin{align}
&\langle f, E_{\opDhat}((\lambda_1,\lambda_2])g\rangle_{L^2(\bbR)^2}  \\
&\ =\lim_{\delta \downarrow 0}\lim_{\varepsilon \downarrow
0}\frac1{2\pi i}\int_{\lambda_1+\delta}^{\lambda_2+\delta}d\lambda
\, \big\langle f, U^{-1}\big[(\opHhat-(\lambda+i\varepsilon))^{-1}-
(\opHhat-(\lambda-i\varepsilon))^{-1}\big]Ug
\big\rangle_{L^2(\bbR)^2}.  \notag
\end{align}
\end{proof}

\subsection{Examples.}\label{s3.4}

We now consider the calculation of quantities discussed in the
previous section for $\Hhat$ and $\Dhat$ for the special case where
$q(x)=q_0\in\bbC$ is constant. In this case we denote $\Hhat$ and $\Dhat$ by
$\Hhat_{q_0}$ and $\Dhat_{q_0}$, etc. But first we consider the case $q_0=0$
and denote $\Hhat$ and $\Dhat$ by $\Hhat_0$ and $\Dhat_0$, etc.

\medskip

\noindent $(i)$ The case $q_0=0$: \\
By direct calculation
for general $\alpha=(\cos(\theta),\sin(\theta))$, $\theta\in [0,2\pi)$, and
$z\in\bbC\backslash\bbR$,
\begin{equation}
\Psi_\pm^{{\Hhat}_{0}}(z,x,\alpha)=
\begin{cases}
a_\pm\begin{pmatrix}1\\\pm i\end{pmatrix}e^{\pm izx}, & \Im (z) >0, \\[15pt]
a_\pm\begin{pmatrix}1\\\mp i\end{pmatrix}e^{\mp izx}, & \Im (z) <0
\end{cases}
\end{equation}
for some $a_\pm \in\bbC$. As noted earlier,
\begin{equation}
\Psi^{\Dhat_0}_\pm (z,x,\beta)=U^{-1}\Psi_\pm^{\Hhat_0}(z,x,\alpha)
\end{equation}
for the corresponding general $\beta=\alpha
U=[(-1+i)/2](e^{-i\theta},e^{i\theta})$, $\theta\in [0,2\pi)$,
where $U$ is defined in \eqref{U}. Explicitly,
\begin{align}
& \Psi^{\Dhat_0}_+ (z,x,\beta)=b_+\begin{pmatrix} 0 \\  1+i\end{pmatrix}
e^{izx}, \quad  \Psi^{\Dhat_0}_- (z,x,\beta)=b_-\begin{pmatrix} 1+i \\
0 \end{pmatrix} e^{-izx}, \quad \Im(z)>0, \\
& \Psi^{\Dhat_0}_+ (z,x,\beta)=b_+ \begin{pmatrix} 1+i \\  0 \end{pmatrix}
e^{-izx}, \quad  \Psi^{\Dhat_0}_- (z,x,\beta)=b_-\begin{pmatrix} 0 \\
1+i \end{pmatrix} e^{izx}, \quad \Im(z)<0
\end{align}
for some $b_\pm \in\bbC$.

In particular, for $\alpha=\alpha_0=(1,0)$ we see that
$\calF^{\Hhat_0}(z,0,\alpha_0) = I_2$ and hence by \eqref{3.15} that
\begin{equation}
\Psi^{\Hhat_0}_\pm(z,0,\alpha_0) =\begin{pmatrix}1\\
\mhat_{0,\pm}(z,\gamma_0)\end{pmatrix}, \quad z\in\bbC\backslash\bbR.
\end{equation}
           From this we conclude that $a_\pm=1$ for $\alpha=\alpha_0$ and that
\begin{equation}
\mhat_{0,\pm}(z,\gamma_0)=m^{\Hhat_{0}}_\pm(z,\alpha_0)
=m^{\Dhat_{0}}_\pm(z,\beta_0)
=\begin{cases} \pm i, & \Im (z) >0,\\ \mp i,
& \Im (z) < 0.
\end{cases}
\end{equation}
As a consequence, the whole-line Weyl--Titchmarsh $M$-function defined in
\eqref{Mhat} is given by
\begin{equation}
\hatt M_0(z,\gamma_0)=\pm(i/2)I_2, \quad \Im (z) \gtrless 0.
\end{equation}
Hence, for $q_0=0$, the spectral measure for $\opHhat_0$, as
described in Theorem~\ref{tWSHhat}, is given by
\begin{equation}
d\hatt \Omega_0(\lambda,\gamma_0)=[1/(2\pi)] I_2\, d\lambda.
\end{equation}

By \eqref{3.12a}, we see that
\begin{equation}
G^{\hatt {H}_0} (z,x,x')= \f{1}{2}
\begin{pmatrix}i&
\mp 1\\
\pm 1 & i \end{pmatrix}e^{\pm iz(x'-x)}, \quad x\lessgtr x', \quad
\Im(z)>0.
\end{equation}
For the corresponding Green's matrix $G^{\hatt{D}_0}(z,x,x')=
U^{-1}G^{\hatt{H}_0}(z,x,x')U$ one obtains
\begin{equation}
G^{\hatt{D}_0}(z,x,x')=
\begin{cases}
i \begin{pmatrix}1&0\\0&0 \end{pmatrix}e^{ iz(x'-x)}, & x<x',\\[15pt]
i \begin{pmatrix}0&0\\0&1 \end{pmatrix}e^{- iz(x'-x)}, & x>x',
\end{cases} \quad \Im(z)>0.
\end{equation}
Similarly,
\begin{equation}
G^{\hatt{H}_0} (z,x,x')= \frac{1}2 \begin{pmatrix}-i&
\mp 1\\
\pm 1 & -i \end{pmatrix}e^{\mp iz(x'-x)},\quad x\lessgtr x', \quad \Im(z)<0,
\end{equation}
and
\begin{equation}
G^{\hatt{D}_0}(z,x,x')=
\begin{cases}
-i \begin{pmatrix}0&0\\0&1 \end{pmatrix}e^{- iz(x'-x)}, & x<x',\\[15pt]
-i \begin{pmatrix}1&0\\0&0 \end{pmatrix}e^{ iz(x'-x)}, & x>x', \
\end{cases} \quad \Im(z)<0.
\end{equation}

The spectra of $\hatt{H}_{0}$ and $\hatt{D}_{0}$ are purely absolutely
continous of uniform multiplicity two and given by
\begin{equation}
\sigma(\hatt{H}_{0})=\sigma (\hatt{D}_{0})=\bbR.
\end{equation}

\medskip

\noindent $(ii)$ The case $q_0\in\bbC\backslash\{0\}$: \\
In considering the case where $q=q_0$ is a nonzero complex constant,
we first define $\hatt S_{q_0}(z)$ to be a function that is analytic with
positive imaginary part on the split plane
\begin{equation}
\hatt\bbP_{q_0}= \bbC\backslash
\{ \lambda\in\bbR \,|\, |\lambda|\ge |q_0| \},
\end{equation}
such that
\begin{equation} \label{3.35}
\hatt S_{q_0}(z) = \sqrt{z^2 - |q_0|^2},
\quad \Im\big(\hatt S_{q_0}(z)\big)>0, \quad z\in\hatt\bbP_{q_0}.
\end{equation}
Thus,
\begin{equation}
\hatt S_{q_0}(\ol z)=- \ol{\hatt S_{q_0}(z)}, \quad z\in \hatt\bbP_{q_0},
\end{equation}
\begin{equation}
\hatt S_{q_0}(\lambda\pm i0)= \lim_{\varepsilon\downarrow 0}
\hatt S_{q_0}(\lambda
\pm i\varepsilon) =
\begin{cases}
\pm\sqrt{\lambda^2-|q_0|^2}, & \lambda \ge |q_0|,\\
\mp\sqrt{\lambda^2-|q_0|^2}, & \lambda \le -|q_0|.
\end{cases}
\end{equation}
(If $q_0=0$, this convention amounts to defining $\sqrt{z^2}=\pm z$ for
$\Im(z) \gtrless 0$.)

A direct calculation shows for general
$\alpha=(\cos(\theta),\sin(\theta))$, $\theta\in [0,2\pi)$, that
\begin{equation}
\Psi^{\Hhat_{q_0}}_\pm(z,x,\alpha) = a_\pm \begin{pmatrix}1\\
\frac{-\Re(q_0)
\pm i \hatt S_{q_0}(z)}{z+\Im(q_0)}
\end{pmatrix}e^{\pm i \hatt S_{q_0}(z) x}, \quad z\in\bbC\backslash\bbR
\end{equation}
for some $a_\pm \in\bbC$. For the corresponding general
$\beta=\alpha U=[(-1+i)/2](e^{-i\theta},e^{i\theta})$, $\theta\in
[0,2\pi)$, we  see by direct calculation that
\begin{equation}
\Psi^{\Dhat_{q_0}}_\pm(z,x,\beta)=
b_\pm\begin{pmatrix}1\\ \frac{i}{q_0}\big[z\pm
\hatt S_{q_0}(z)\big]\end{pmatrix} e^{\pm i \hatt S_{q_0}(z)x}, \quad
z\in\bbC\backslash\bbR
\end{equation}
for some $b_\pm \in\bbC$, and alternatively that
\begin{align}
\begin{split}
\Psi^{\Dhat_{q_0}}_\pm(z,x,\beta)&=U^{-1}
\Psi^{\Hhat_{q_0}}_\pm(z,x,\alpha)\\
&=\frac{a_\pm(1+i)}{2(z+\Im(q_0))}
\begin{pmatrix}
-z+iq_0 \pm \hatt S_{q_0}(z)\\
-z-i\ol q_0 \mp \hatt S_{q_0}(z)
\end{pmatrix}e^{\pm iS_{q_0}(z)x}, \quad z\in\bbC\backslash\bbR.
\end{split}
\end{align}
In particular, for $\alpha=\alpha_0=(1,0)$ we see that
$\calF^{\hatt{H}_{q_0}}(z,0,\alpha_0) = I_2$ and hence by
\eqref{3.15} that
\begin{equation}
\Psi^{\Hhat_{q_0}}_\pm(z,0,\alpha_0) =\begin{pmatrix}1\\
\mhat_{q_0,\pm}(z,\gamma_0)\end{pmatrix} =a_\pm \begin{pmatrix}1\\
\frac{-\Re(q_0) \pm i \hatt S_{q_0}(z)}{z+\Im(q_0)} \end{pmatrix},
\quad z\in\bbC\backslash\bbR.
\end{equation}
      From this we conclude that $a_\pm=1$ for $\alpha=\alpha_0$ and that
\begin{equation}
\mhat_{q_0,\pm}(z,\gamma_0)=m^{\Hhat_{q_0}}_\pm(z,\alpha_0)
=m^{\Dhat_{q_0}}_\pm(z,\beta_0)=\frac{-\Re(q_0)
\pm i \hatt S_{q_0}(z)}{z+\Im(q_0)}, \quad z\in\bbC\backslash\bbR.
\end{equation}

For $\Re(q_0)\ne 0$, we note that $\mhat_{q_0,\pm}(z,\gamma_0)$ is
analytic for $z\in (-|q_0|, |q_0|)$ with the possible exception of
$z=-\Im (q_0)$. In this case, the $z$-wave functions of
$\Hhat_{q_0}$ corresponding to $z=-\Im (q_0)$ are given by
\begin{equation}
\Psi^{\Hhat_{q_0}}(x)= \psi^{\Hhat_{q_0}}_1(0)\begin{pmatrix}e^{x\Re(q_0) } \\
2\left(\frac{\Im (q_0)}{\Re(q_0)}\right)\sinh (x \, \Re(q_0)
)\end{pmatrix} +\psi^{\Hhat_{q_0}}_2(0)\begin{pmatrix}0\\ e^{-x\Re(
q_0)}
\end{pmatrix}.
\end{equation}
As a consequence, we see that while $z=-\Im (q_0)$ is not an
eigenvalue for $\opHhat_{q_0}$, it is an eigenvalue for
$\opHhat_{q_0,\pm}(\alpha_0)$ corresponding to a simple pole for
$\mhat_{q_0,\pm}(z,\gamma_0)$ for $\Re(q_0)\gtrless 0$. We also note
that $z=-\Im (q_0)$ corresponds to a removable singularity for
$\mhat_{q_0,\mp}(z,\gamma_0)$ for $\Re(q_0)\gtrless 0$.

However, for $\Re(q_0)= 0$, $z=-\Im (q_0)$ corresponds to an
endpoint of the spectral gap $(-|q_0|, |q_0|)$ and  the $z$-wave
functions of $\Hhat_{q_0}$ are given by
\begin{equation}
\Psi^{\Hhat_{q_0}}(x)= \psi^{\Hhat_{q_0}}_1(0)\begin{pmatrix}1\\2x \,
\Im (q_0) \end{pmatrix} + \psi^{\Hhat_{q_0}}_2(0)\begin{pmatrix}0\\1
\end{pmatrix}.
\end{equation}
In this case, $\Psi^{\Hhat_{q_0}}(x)$ is neither in $L^2(\mathbb
R)^2$ nor in $L^2([0,\pm\infty))^2$; hence $z=-\Im (q_0)$ is not an
eigenvalue for $\opHhat_{q_0}$, or for
$\opHhat_{q_0,\pm}(\alpha_0)$.

As a consequence, the whole-line Weyl--Titchmarsh $M$-function defined in
\eqref{Mhat} is now given by
\begin{equation}
\Mhat_{q_0}(z,\alpha_0)=\frac{i}{2\hatt
S_{q_0}(z)}\begin{pmatrix}z+\Im(q_0)& -\Re(q_0)\\ -\Re(q_0)&
z-\Im(q_0)\end{pmatrix}, \quad z\in\bbC\backslash\bbR.
\end{equation}
Hence, the spectral measure for $\opHhat_{q_0}$, as described in
Theorem~\ref{tWSHhat} by $d\hatt\Omega(\lambda,\gamma)$  in \eqref{3.42},
is determined by $\lim_{\varepsilon\downarrow 0}\Im\big(\hatt
M_{q_0}(\lambda +i\varepsilon,\gamma_0 )\big)$ and, in light of
\eqref{3.35}, found to be
\begin{equation}
d\hatt\Omega_{q_0}(\lambda,\gamma_0)=d\lambda
\begin{cases}
\begin{pmatrix} 0 &0 \\ 0 & 0 \end{pmatrix}, & |\lambda|<|q_0|\\
\dfrac{\pm 1}{2\pi \hatt S_{q_0}(\lambda+i0)}
\begin{pmatrix}
\lambda + \Im(q_0) & - \Re(q_0)\\
-\Re(q_0) & \lambda - \Im(q_0)
\end{pmatrix}, & \pm\lambda \ge |q_0|.
\end{cases}
\end{equation}

By \eqref{3.12a} we see that
\begin{align}
&G^{\hatt{H}_{q_0}}(z,x,x')  \\
&\quad =\frac{i}{2\hatt S_{q_0}(z)}\begin{pmatrix}z+\Im(q_0)& -\Re(q_0)\pm
i\hatt S_{q_0}(z)\\
-\Re(q_0)\mp i\hatt S_{q_0}(z)&
z-\Im(q_0)\end{pmatrix}e^{\pm i\hatt S_{q_0}(z)(x'-x)},
\quad x\lessgtr x', \no \\
& \hspace*{10.95cm} z\in\bbC\backslash\bbR. \no
\end{align}
Using  either \eqref{3.16b} or the fact that
$G^{\hatt{D}_{q_0}}(z,x,x')=U^{-1}G^{\hatt{H}_{q_0}} (z,x,x')U$, we obtain
\begin{align}
& G^{\hatt{D}_{q_0}} (z,x,x')=\frac{1}{2\hatt S_{q_0}(z)}
\begin{pmatrix}
i\big[z\pm \hatt S_{q_0}(z)\big] & q_0 \\
-\ol{q_0} & i\big[z\mp \hatt S_{q_0}(z)\big]
\end{pmatrix}e^{\pm i \hatt S_{q_0}(z)(x'-x)},\quad x\lessgtr x', \\
& \hspace*{10.85cm} z\in\bbC\backslash\bbR. \no
\end{align}
The spectra of $\hatt{H}_{q_0}$ and $\hatt{D}_{q_0}$ are purely absolutely
continous of uniform multiplicity two and given by
\begin{equation}
\sigma(\hatt{H}_{q_0})=\sigma (\hatt{D}_{q_0})=(-\infty,-|q_0|]\cup
[|q_0|,\infty).
\end{equation}

\section{Non-self-adjoint Dirac and Hamiltonian Systems}  \lb{s4}

In this section, we focus upon $J$-self-adjoint realizations for $D$
and its unitarily equivalent $H$, specifically, the operator
$\opDwti$, defined by \eqref{DO} corresponding to the special
case of \eqref{DE} given by \eqref{DENSA}, and the operator
$\opHwti$ defined by \eqref{HO} associated with the special
case of \eqref{HE} given by \eqref{HENSA}. Some spectral theory for
the non-self-adjoint operator $\opDwti$, and therefore for its
unitary equivalent $\opHwti$, has been developed in \cite{CGHL04}.
However, it remains incomplete by comparison with their self-adjoint
counterparts $\opDhat$ and $\opHhat$ as described in the previous
section.

\subsection{Weyl--Titchmarsh coefficients.}

We now turn to the development in the non-self-adjoint setting of
the analog for the Weyl--Titchmarsh coefficient defined and
discussed in Section~\ref{sWTCSA}. This subsection details (and
partially corrects) Remark 5.6 in \cite{CGHL04} which anticipated
the introduction of half-line Weyl--Titchmarsh $m$-functions
associated with $\wti \scrD$. We note that a general
Weyl--Titchmarsh--Sims theory for singular non-self-adjoint Hamiltonian
systems has recently been developed in \cite{BEP03} (see also 
\cite{BM03} for additional spectral results and further references).
However, while the general case considered in \cite{BEP03} requires
certain restrictions on the complex spectral parameter $z$ when introducing
a Weyl--Titchmarsh coefficient $m(z)$, the very special structure of $\wti
\scrD$ permits us to introduce a Weyl--Titchmarsh coefficient on the
resolvent set $\rho(\opDwti)=\rho(\opHwti)$ in this section. We also
emphasize that the Weyl--Titchmarsh $m$-coefficient was
first introduced for a class of $J$-self-adjoint Dirac-type operators with
bounded coefficients (and for the complex spectral parameter restricted to
a half-plane) in \cite{Sa90} (see also \cite{GKS98}, \cite{Sa05} and the
literature therein).

\begin{hypothesis}
Throughout this section, we assume that the resolvent set
$\rho(\opDwti)$ of $\opDwti$ $($and hence that of $\opHwti$$)$ is
nonempty.
\end{hypothesis}
To begin, we note the fundamental result
established in \cite[Theorem 5.4]{CGHL04} which states that
\begin{equation}\label{wtiNulla}
\dim \big(\DwtiN(z,\pm\infty)\big)=1, \quad z\in\rho(\opDwti),
\end{equation}
and hence by the unitary equivalence given in \eqref{HE},
\begin{equation}
\dim\big(\HwtiN(z,\pm\infty)\big)=1,\quad z\in\rho(\opHwti)=\rho(\opDwti).
\end{equation}
In particular, one has the following result.

\begin{theorem} [\cite{CGHL04}] \lb{t4.1.1}
The operator $\opDwti$, maximally defined in \eqref{DO}, is
$J$-self-adjoint since
\begin{equation}
\calJ \opDwti \calJ=\opDwti^*,  \lb{4.4}
\end{equation}
where $\calJ$ is defined in \eqref{calJ}, and the operator $\opHwti$,
maximally defined in \eqref{HO}, is
$J$-self-adjoint since
\begin{equation}
\wti\calJ \opHwti \wti\calJ=\opHwti^*,   \lb{4.3}
\end{equation}
where $\wti\calJ$ denotes the conjugate linear involution
\begin{equation}
\wti\calJ= i\calC I_2.  \lb{4.5}
\end{equation}
Moreover, $\opHwti$ and $\opDwti$ are unitarily equivalent, i.e.
\begin{equation}
\opHwti=U\opDwti U^{-1}.  \lb{4.6a}
\end{equation}
\end{theorem}
\begin{proof}
Equation \eqref{4.3} has been proven in \cite{CGHL04}. The rest follows
from the unitary equivalence \eqref{HE} via the constant unitary matrix $U$.
\end{proof}
As in the self-adjoint setting, one defines the half-line operator
$\opDwti_\pm(\beta)$ in association with the differential expression
$\Dwti$ found in \eqref{DENSA} by
\begin{align}
\opDwti_\pm(\beta)&=i\sigma_3\dfrac{d}{dx} + \wti Q,  \no \\
\dom(\opDwti_\pm(\beta))&= \big\{ \Psi\in
L^2([0,\pm\infty))^2\,\big|\, \Psi\in\AC_\loc([0,\pm\infty))^2, \\
& \hspace*{2cm} \beta\Psi(0)=0,\ \Dwti\Psi\in
L^2([0,\pm\infty))^2\big\}, \no
\end{align}
where
\begin{equation}
\beta=[(-1+i)/2] (e^{-i\theta},e^{i\theta}), \quad
\theta\in [0,2\pi),
\end{equation}
and where $\beta\Psi(0)=0$ represents a $J$-self-adjoint boundary
condition for $\opDwti_\pm(\beta)$ using the conjugation $\calJ$.
One also defines the half-line operator $\opHwti_\pm(\alpha)$ in
association with the differential expression $\Hwti$ found in
\eqref{HENSA} by
\begin{align}
\opHwti_\pm(\alpha)&=-\sigma_4\dfrac{d}{dx} + \wti B,  \no \\
\dom(\opHwti_\pm(\alpha))&=\big\{ \Psi\in
L^2([0,\pm\infty))^2 \,\big|\, \Psi\in\AC_\loc([0,\pm\infty))^2, \\
&\hspace*{2cm} \alpha\Psi(0)=0,\ \Hwti\Psi\in
L^2([0,\pm\infty))^2\big\}, \no
\end{align}
where $\alpha=\beta U^{-1}=(\cos(\theta),\sin(\theta))$, $\theta\in
[0,2\pi)$, and where $\alpha\Psi(0)=0$ represents a $J$-self-adjoint
boundary condition for $\opHwti_\pm(\alpha)$ using the conjugation
$\wti\calJ$.  In fact, $\opDwti_\pm(\beta)$ and $\opHwti_\pm(\alpha)$ are
also
$J$-self-adjoint,
\begin{equation}
\calJ \opDwti_\pm(\beta) \calJ = \opDwti_\pm(\beta)^*, \quad
\wti\calJ \opHwti_\pm(\alpha) \wti\calJ = \opHwti_\pm(\alpha)^*. \lb{4.10}
\end{equation}

To prove \eqref{4.10} one first notes that apart from the
$J$-self-adjoint boundary condition imposed at $x=0$,
$\opDwti_\pm(\beta)$ and $\opHwti_\pm(\alpha)$ are maximally defined
and one only needs to check the corresponding $L^2([0,\pm\infty))^2$
condition in a neighborhood of $\pm\infty$. But the latter
immediately follows from \eqref{4.4} and \eqref{4.3}. As in the
self-adjoint context (cf.\ \eqref{3.6a}) one infers that
\begin{equation}
\opHwti_\pm(\alpha)=U \opDwti_\pm(\beta) U^{-1}, \quad \beta=\alpha U,
\end{equation}
holds in addition to \eqref{4.6a}.

\begin{remark} \lb{r4.1.2}
Also exploited in \cite[Lemma 5.2, Theorem 5.4]{CGHL04} is a feature
that distinguishes
$\opDwti$ from $\opDhat$: The bijection $\calK$ acting upon
$\AC_{\loc}(\bbR)^2$ and described by
\begin{equation}\label{calK}
\calK = \sigma_4\calC, \quad \calK^2=-I_2,
\end{equation}
where $\calC$ is the conjugation operator acting on $\bbC^2$ defined
in \eqref{Conj}, maps $z$-wave functions of $\Dwti$ to $\bar z$-wave
functions of $\Dwti$. By this, we mean that
\begin{equation}\label{wavefcnmap}
(\Dwti\Psi)(z,x)=z\Psi(z,x) \, \text{ if and only if } \,
(\Dwti\calK\Psi)(z,x)=\bar z\calK\Psi(z,x), \quad z\in\bbC, \; x\in\bbR.
\end{equation}
By contrast, $\calK$ fails to map $z$-wave functions to $\bar
z$-wave functions of $\Dhat$. Distinguishing $\Dhat$ from $\Dwti$
is the fact that rather than $\calK$, it is the operator  $\calJ$,
defined in \eqref{calJ} and acting as a bijection on
$\AC_{\loc}(\bbR)^2$, that serves to map $z$-wave functions to $\bar
z$-wave functions of $\Dhat$:
\begin{equation}
(\Dhat \Psi)(z,x)=z\Psi(z,x) \, \text{ if and only if } \,
(\Dhat\calJ\Psi)(z,x)=\bar z\calJ \Psi(z,x), \quad z\in\bbC, \; x\in\bbR.
\end{equation}
\end{remark}

As before, we  introduce the fundamental system of solutions of
$\Hwti\Psi=z\Psi$ by
\begin{equation}
\HwtiTheta(z,\cdot,\alpha),\, \HwtiPhi(z,\cdot,\alpha) \in \AC_\loc(\bbR)^2,
\quad z\in\bbC,
\end{equation}
and the
matrix-valued function $\HwtiF(z,\cdot,\alpha)$  given by
\begin{equation} \label{FSHwti}
\HwtiF(z,x,\alpha)=\big(\HwtiTheta(z,x,\alpha)\ \ \HwtiPhi(z,x,\alpha)\big),
\end{equation}
where for $\theta\in [0,2\pi)$,
\begin{equation}\label{FSHwtib}
\HwtiTheta(z,0,\alpha)=(\cos(\theta), \sin(\theta))^\top ,\quad
\HwtiPhi(z,0,\alpha)= (-\sin(\theta), \cos(\theta))^\top.
\end{equation}
We also introduce the related fundamental system of $z$-wave
functions of $\Dwti$ given by
$\DwtiTheta(z,\cdot,\beta)$ and  $\DwtiPhi(z,\cdot,\beta)$, as well as
the matrix-valued function
\begin{equation}
\DwtiF(z,x,\beta)=\big(\DwtiTheta(z,x,\beta)\ \ \DwtiPhi(z,x,\beta)
\big)=U^{-1}\HwtiF(z,x,\alpha),
\end{equation}
where $\beta=\alpha U$. Thus, for $\theta\in [0,2\pi)$,
\begin{equation}\label{FSDwti}
\DhatTheta(z,0,\beta)=[(1+i)/2] (e^{i\theta},e^{-i\theta})^\top ,\quad
\DhatPhi(z,0,\beta)=[(1-i)/2] (e^{i\theta},-e^{-i\theta})^\top.
\end{equation}

Analogous to the self-adjoint setting, a Weyl-Titchmarsh coefficient
can be defined for values of $z\in\bbC$ that lie in the compliment
of the combined spectrum for $\opDwti$ and $\opDwti_\pm(\beta)$. The
fact that $\DwtiPsi_\pm(z,\cdot,\beta)$ form a basis for the
$z$-wave functions of $\Dwti$ implies that
$\DwtiPsi_\pm(z,\cdot,\beta)$ are not scalar multiples of
$\DwtiPhi(z,\cdot,\beta)$ for $z\in
\rho(\opDwti)\cap\rho(\opDwti_\pm(\beta))$ and hence similarly
that $\HwtiPsi_\pm(z,\cdot,\alpha)$ are not scalar multiples of
$\HwtiPhi(z,\cdot,\alpha)$ for $z\in\
\rho(\opDwti)\cap\rho(\opDwti_\pm(\beta))=
\rho(\opHwti)\cap\rho(\opHwti_\pm(\alpha))$ .  This being the
case, $\HwtiPsi_\pm$ and $\DwtiPsi_\pm$ have the unique
representations given by
\begin{align}
\HwtiPsi_\pm (z,\cdot,\alpha) &=\HwtiTheta(z,\cdot,\alpha) +
\Hwtim_\pm (z,\alpha)\HwtiPhi(z,\cdot,\alpha)\in
L^2([0,\pm\infty))^2, \label{WSHwti} \\
\DwtiPsi_\pm (z,\cdot,\beta) &=\DwtiTheta(z,\cdot,\beta) +
\Dwtim_\pm
(z,\beta)\DwtiPhi(z,\cdot,\beta)\in L^2([0,\pm\infty))^2. \label{WSDwti}
\end{align}
In complete analogy to Lemma~\ref{l3.1}, and by completely analogous
proof, one obtains the following result:
\begin{lemma}
Let $\alpha=(\cos(\theta),\sin(\theta))$, $\theta\in [0,2\pi)$, and
let $\beta = \alpha U$ with $U$ defined in \eqref{U}. Let
$\HwtiTheta$, $\HwtiPhi$ represent the fundamental system of
solutions of the Hamiltonian system $\Hwti\Psi=z\Psi$ satisfying
\eqref{FSHwtib}, and let $\DwtiTheta$, $\DwtiPhi$ represent the
fundamental system of solutions of the Dirac system
$\Dwti\Psi=z\Psi$ satisfying \eqref{FSDwti}. Then, for $z\in
\rho(\opDwti)\cap\rho(\opDwti_\pm(\beta))=
\rho(\opHwti)\cap\rho(\opHwti_\pm(\alpha))$, with $\HwtiPsi_\pm$
defined in \eqref{WSHwti} and with $\DwtiPsi_\pm$ defined in
\eqref{WSDwti}, one infers that
\begin{equation}
\DwtiPsi_\pm(z,x,\beta)=U^{-1}\HwtiPsi_\pm(z,x,\alpha), \quad
x\in\bbR,
\end{equation}
and in particular, that
\begin{equation}\label{4.14}
\Dwtim_\pm(z,\beta)=\Hwtim_\pm(z,\alpha).
\end{equation}
\end{lemma}

\begin{remark}  \lb{r4.4}
As in Remark~\ref{r3.2}, we denote
$\mwti_{\pm}(z,\gamma)=\Dwtim_\pm(z,\beta) =\Hwtim_{\pm}(z,\alpha)$, where
$\gamma$ represents $\alpha$ in the context of $\wti\scrH$ and
$\beta=\alpha U$ in the context of $\wti\scrD$.
\end{remark}

It is well-known that $\ol{\mhat_\pm(z,\gamma)}=\mhat_\pm(\ol z,\gamma)$,
$z\in\bbC\backslash\bbR$, for the self-adjoint differential expression
$\Dhat$. By contrast, this is not the case in the non-self-adjoint setting
for $\Dwti$ as seen in the next result.
\begin{lemma}\label{l4.4}
Let $\alpha=\alpha(\theta)=(\cos(\theta),\sin(\theta))$,
$\beta=\beta(\theta)=[(-1+i)/2] (e^{-i\theta},e^{i\theta})$, and
$\gamma(\theta)$ defined as in Remark \ref{r4.4},
$\theta\in [0,2\pi)$, and $z\in
\rho(\opDwti)\cap\rho(\opDwti_\pm(\beta))$. Then,
\begin{align}
\ol{\mwti_\pm(z,\gamma(\theta))}&=-[\mwti_\pm(\ol z,
\gamma(\theta))]^{-1},   \label{4.16}\\
\ol{\Dwtim_\pm(z,\beta(\theta))}&=\Dwtim_\pm(\ol z,
\beta(\theta-\pi/2)),  \label{4.17}\\
\ol{\Hwtim_\pm(z,\alpha(\theta))}&=\Hwtim_\pm(\ol z,
\alpha(\theta-\pi/2)).  \label{4.18}
\end{align}
\end{lemma}

\begin{proof}
As defined in \eqref{calK}, $\calK$ is an isometry on
$L^2([0,\pm\infty))^2$ which by \eqref{wavefcnmap} maps $z$-wave
functions to $\ol z$-wave functions of $\Dwti$. As a consequence of
\eqref{wtiNulla},
$\calK\DwtiPsi_\pm(z,\cdot,\beta)=c_\pm\DwtiPsi_\pm(\ol
z,\cdot,\beta)$ for some $c_\pm\in\bbC$. Moreover, for the
fundamental system defined in \eqref{FSDwti},
$\calK\DwtiTheta(z,0,\beta)=-\DwtiPhi(\ol z,0,\beta),\
\calK\DwtiPhi(z,0,\beta)=\DwtiTheta(\ol z,0,\beta)$, and hence
\begin{equation}
\calK\DwtiTheta(z,x,\beta)=-\DwtiPhi(\ol z,x,\beta),\quad
\calK\DwtiPhi(z,x,\beta)=\DwtiTheta(\ol z,x,\beta), \quad x\in\bbR.
\end{equation}
Given the unique representations provided by \eqref{WSDwti}, we
obtain for $\theta_1, \theta_2\in [0,2\pi)$, that
\begin{align}
&c_\pm(\DwtiTheta(\ol z,x,\beta(\theta_1)) + \Dwtim_\pm
(\ol z,\beta(\theta_1))\DwtiPhi(\ol z,x,\beta(\theta_1))) \no \\
&\quad=\calK(\DwtiTheta(z,x,\beta(\theta_2)) + \Dwtim_\pm
(z,\beta(\theta_2))\DwtiPhi(z,x,\beta(\theta_2)))\notag\\
&\quad =-\DwtiPhi(\ol z,x,\beta(\theta_2))+\ol{\Dwtim_\pm
(z,\beta(\theta_2))}\DwtiTheta(\ol z,x,\beta(\theta_2)).
\end{align}
From this we see that
\begin{equation}
c_\pm\begin{pmatrix}1\\ \Dwtim_\pm (\ol z,\beta(\theta_1))\end{pmatrix}=
\DwtiF(\ol z,x,\beta(\theta_1))^{-1}\DwtiF(\ol z,x,\beta(\theta_2))
\begin{pmatrix}\ol{\Dwtim_\pm (z,\beta(\theta_2))}\\ -1
\end{pmatrix}.
\end{equation}
Then, \eqref{4.16} follows from $\DwtiF(\ol
z,x,\beta(\theta))^{-1}\DwtiF(\ol z,x,\beta(\theta))=I_2$;
\eqref{4.17} follows from $\DwtiF(\ol
z,0,\beta(\theta-\pi/2))^{-1}\DwtiF(\ol z,0,\beta(\theta))=J$; and
\eqref{4.18} follows from \eqref{4.14}.
\end{proof}

\begin{remark}
The same argument used to prove \eqref{4.16} can be used in the
self-adjoint context to prove that
\begin{equation}
\ol{\mhat_\pm(z,\gamma)}=\mhat_\pm (\ol z,\gamma), \quad
z\in\bbC\backslash\bbR,
\end{equation}
where the operator $\calJ$ which maps $z$-wave functions to
$\ol z$-wave functions of $\Dhat$ is used rather than $\calK$.
\end{remark}

\subsection{Green's matrices.}

Given the existence of the half-line Weyl--Titchmarsh solutions
$\DwtiPsi_\pm(z,\cdot,\beta)\in L^2([0,\pm\infty))^2$ for
$\Dwti\Psi=z\Psi$, and the corresponding solutions
$\HwtiPsi_\pm(z,\cdot,\alpha)=U\DwtiPsi_\pm(z,\cdot,\beta)\in
L^2([0,\pm\infty))^2$ of $\Hwti\Psi=z\Psi$ for
$z\in\rho(\opHwti)=\rho(\opDwti)$, we now obtain, as special cases of
Lemmas \ref{l2.1} and \ref{l2.3}, descriptions for the Green's matrices
corresponding to the whole-line operators $\opDwti$ and
$\opHwti$:

\begin{lemma}\label{l4.2}
With $z\in\rho(\opHwti)$, let $\HwtiPsi_\pm(z,\cdot,\alpha)\in
L^2([0,\pm\infty))^{2}$ represent a basis of solutions for the
Hamiltonian system $\Hwti\Psi=z\Psi$. Then, the whole-line Green's
matrix associated with the operator $\opHwti$ is given by
\begin{equation}\label{4.7}
\GopHwti(z,x,x')= K(z,\alpha)\HwtiPsi_\mp(z,x,\alpha)
\HwtiPsi_\pm(z,x',\alpha)^\top, \quad x\lessgtr x',
\end{equation}
where
\begin{align}
K(z,\alpha)&=\big[W(\HwtiPsi_+(z,x,\alpha),
\HwtiPsi_-(z,x,\alpha))\big]^{-1}
=\big[W(\DwtiPsi_+(z,x,\beta),\DwtiPsi_-(z,x,\beta))\big]^{-1} \\
&= [\wti m_-(z,\gamma) - \wti m_+(z,\gamma)]^{-1}
\end{align}
is constant with respect to $x\in\bbR$.
\end{lemma}

\begin{lemma}\label{l4.1}
With $z\in\rho(\opDwti)$, let $ \DwtiPsi_{\pm}(z,\cdot,\beta)\in
L^2([0,\pm\infty))^{2}$ represent a basis of z-wave functions of
$\Dwti$. Then, the whole-line Green's matrix  associated with the
operator $\opDwti$ is given by
\begin{equation}\label{4.6}
\GopDwti(z,x,x')=C(z,\beta) \DwtiPsi_\mp(z,x,\beta)
\DwtiPsi_\pm(z,x',\beta)^\top \sigma_1,\quad x\lessgtr x'
\end{equation}
where $\sigma_1$ is given in \eqref{calJ} and where
\begin{equation}
C(z,\beta)=-i\big[W(\DwtiPsi_+(z,x,\beta),\DwtiPsi_-(z,x,\beta))\big]^{-1}
\end{equation}
is constant with respect to $x\in\bbR$.
\end{lemma}

With $\mwti_\pm(\cdot,\gamma)$ defined in the non-self-adjoint settings of
$\Dwti$ and the unitarily equivalent $\Hwti$, we define
$\Gammawti(\cdot,\gamma)$  by substituting $\mwti_\pm(\cdot,\gamma)$ for
its corresponding $\mhat_\pm(\cdot,\gamma)$ in the definition of
$\Gammahat(\cdot,\gamma)$ given in \eqref{Gammahat}. That is,
\begin{align}
& \Gammawti(z,\gamma)=\begin{pmatrix}
\frac1{\mwti_-(z,\gamma) - \mwti_+(z,\gamma)}&
\frac{\mwti_-(z,\gamma)}{\mwti_-(z,\gamma) - \mwti_+(z,\gamma)}\\[6pt]
\frac{ \mwti_+(z,\gamma)}{\mwti_-(z,\gamma) - \mwti_+(z,\gamma)}&
\frac{\mwti_-(z,\gamma)\mwti_+(z,\gamma)}{\mwti_-(z,\gamma)
- \mwti_+(z,\gamma)}
\end{pmatrix}, \label{Gammawti} \\
& \hspace*{1.05cm} z\in\rho(\opDwti)\cap\rho(\opDwti_\pm(\beta))=
\rho(\opHwti)\cap\rho(\opHwti_\pm(\alpha)).  \no
\end{align}
With this definition one obtains alternative
expressions for the Green's matrices $\GopDwti(z,x,x')$ and
$\GopHwti(z,x,x')$ given in Lemmas~\ref{l4.1} and \ref{l4.2}  which
are analogous to those given for $\GopDhat(z,x,x')$ and
$\GopHhat(z,x,x')$  in Lemmas~\ref{l3.3} and \ref{l3.4}, that is, for
$z\in \rho(\opDwti)\cap\rho(\opDwti_\pm(\beta))=
\rho(\opHwti)\cap\rho(\opHwti_\pm(\alpha))$,
\begin{align}
\GopHwti(z,x,x')&=\begin{cases}
\HwtiF(z,x,\alpha)\Gammawti(z,\gamma)^\top \HwtiF(z,x',\alpha)^\top ,
& x<x',\\
\HwtiF(z,x,\alpha)\Gammawti(z,\gamma)\HwtiF(z,x',\alpha)^\top , &x>x',
\end{cases}  \lb{4.34} \\
\GopDwti(z,x,x')&=\begin{cases}
-i \DwtiF(z,x,\beta)\Gammawti(z,\gamma)^\top \DwtiF(z,x',\beta)^\top ,
& x<x',\\ -i \DwtiF(z,x,\beta)\Gammawti(z,\gamma) \DwtiF(z,x',\beta)^\top ,
& x>x'.
\end{cases}  \lb{4.35}
\end{align}

\begin{lemma}  \lb{l4.10}
Let $\alpha=(\cos(\theta),\sin(\theta))$, $\theta\in [0,2\pi)$, and let
$\beta = \alpha U$ with $U$ defined in \eqref{U}. Then,
\begin{equation}
\wti m_{\pm}(\cdot,\gamma) \, \text{ are analytic on } \,
\rho(\wti H_{\pm}(\alpha))=\rho(\wti D_{\pm}(\beta)).
\end{equation}
\end{lemma}
\begin{proof}
In analogy to \eqref{4.6}, the half-line Green's matrix of $\wti
H_+ (\alpha)$ is of the form
\begin{equation}\label{4.37}
G^{\wti H_+(\alpha)}(z,x,x') = \HwtiPhi(z,x,\alpha)
\HwtiPsi_+(z,x',\alpha)^\top, \quad 0 \leq x < x', \; z\in
\rho(\wti H_+ (\alpha)).
\end{equation}
Writing
\begin{equation}
\HwtiPhi=(\phi^{\wti \scrH}_1 \; \phi^{\wti \scrH}_2)^\top, \quad
\HwtiTheta=(\theta^{\wti \scrH}_1 \; \theta^{\wti \scrH}_2)^\top, \quad
\HwtiPsi_+=(\psi^{\wti \scrH}_{+,1} \; \psi^{\wti \scrH}_{+,2})^\top,
\lb{4.41}
\end{equation}
we next pick $z_0\in\rho(\wti H_{+}(\alpha))$ and choose $f, g \in
C_0^\infty((0,\infty))$ such that
\begin{equation}
\supp (f) \subseteq [a,b], \quad \supp (g) \subseteq [c,d], \quad
0<a<b<c<d,
\end{equation}
and
\begin{equation}
\int_a^b dx\,
\ol{f(x)} \phi^{\wti \scrH}_1(z_0,x,\alpha)\neq 0, \quad
\int_c^d dx' \, \phi^{\wti \scrH}_1(z_0,x',\alpha) g(x')\neq 0.
\end{equation}
Since for fixed $x\in\bbR$, $\phi^{\wti \scrH}_1(z,x,\alpha)$ and
$\theta^{\wti
\scrH}_1(z,x,\alpha)$ are entire with respect to $z$ and for fixed
$z\in\bbC$ locally absolutely continuous in $x\in\bbR$, continuity with
respect to $z$ then yields
\begin{equation}
\int_a^b dx\,
\ol{f(x)} \phi^{\wti \scrH}_1(z,x,\alpha)\neq 0, \quad
\int_c^d dx' \, \phi^{\wti \scrH}_1(z,x',\alpha) g(x')\neq 0, \quad
z\in \calU(z_0),
\end{equation}
where $\calU(z_0)\subset \rho(\wti H_{+}(\alpha))$ is a sufficiently
small open neighborhood of $z_0$.

Next, one computes for $z\in \calU(z_0)$,
\begin{align}
&\left\langle\begin{pmatrix} f \\ 0 \end{pmatrix}, (\wti
H_{+}(\alpha)-z)^{-1}
\begin{pmatrix} g \\ 0 \end{pmatrix}\right\rangle_{L^2(\bbR)^2}  \lb{4.38}
\\ & \quad = \int_a^b dx\,
\ol{f(x)} \phi^{\wti \scrH}_1(z,x,\alpha)\int_c^d dx' \,
\big[\theta^{\wti \scrH}_1(z,x',\alpha) + \wti m_+(z,\gamma) \phi^{\wti
\scrH}_1(z,x',\alpha)\big] g(x').  \no
\end{align}
Since the left-hand side of \eqref{4.38} is analytic on
$\calU(z_0)$, one concludes that $\wti m_+(\cdot,\gamma)$ is analytic
on $\calU(z_0)$.  Since $z_0\in\rho(\wti H_{+}(\alpha))$ was
arbitrary, $\wti m_+ (\cdot,\gamma)$ (and hence
$G^{\wti H_+(\alpha)}(\cdot,x,x')$, $0\leq x<x'$) is analytic on
$\rho(\wti H_+ (\alpha))$. The analogous proof applies
to $\wti m_- (\cdot,\gamma)$ (and $G^{\wti H_-(\alpha)}(\cdot,x,x')$,
$0\leq x'<x$).
\end{proof}

The argument in the proof of Lemma \ref{l4.10} is a simple variant of the
proof of Lemma 9.1 in \cite{We87} in the context of self-adjoint
higher-order matrix-valued differential operators adapted to the present
case of non-self-adjoint Dirac-type operators. This strategy of
proof markedly differs from the usual approach in the self-adjoint case
which is based on uniform convergence of sequences of Weyl--Titchmarsh
functions lying in nesting Weyl circles. The latter approach generally
fails in the non-self-adjoint context.

A proof analogous to that of Lemma \ref{l4.10} also applies to the
full-line operators $\wti H$ and $\wti D$ and we turn to that next.

\begin{lemma}  \lb{l4.11}
Let $\alpha=(\cos(\theta),\sin(\theta))$, $\theta\in [0,2\pi)$, and let
$\beta = \alpha U$ with $U$ defined in \eqref{U}. Then,
\begin{equation}
\wti \Gamma(\cdot,\gamma) \, \text{ is analytic on } \,
\rho(\wti H)=\rho(\wti D).  \lb{4.45}
\end{equation}
Thus, for $x, x' \in\bbR$, $x\neq x'$, the Green's matrices
$\GopHwti(z,x,x')$ and $\GopDwti(z,x,x')$ are analytic on $\rho(\wti
H)=\rho(\wti D)$ and hence \eqref{Gammawti}--\eqref{4.35} extend to
$\rho(\wti H)=\rho(\wti D)$.
\end{lemma}
\begin{proof}
For simplicity we only consider $\GopHwti(z,x,x')$ for $x<x'$. The case
$x>x'$ and the corresponding results for $\GopDwti(z,x,x')$ follow in an
analogous manner.

Recalling our notation in \eqref{4.41} and suppressing $\alpha$, $\beta$,
and $\gamma$ for simplicity, we start by noting that
\eqref{4.34} yields for the
$(1,1)$-element of
$\GopHwti(z,x,x')$,
\begin{align}
\GopHwti_{1,1}(z,x,x')&= \f{1}{\wti m_-(z)-\wti
m_+(z)}\big[\theta_1^{\wti\scrH}(z,x) \theta_1^{\wti\scrH}(z,x')+\wti
m_+(z) \theta_1^{\wti\scrH}(z,x)\phi_1^{\wti\scrH}(z,x')   \no \\
& \quad\;\;\, + \wti m_-(z)
\phi_1^{\wti\scrH}(z,x)\theta_1^{\wti\scrH}(z,x') + \wti m_-(z)\wti m_+(z)
\phi_1^{\wti\scrH}(z,x)\phi_1^{\wti\scrH}(z,x')\big] \no \\
&= \begin{pmatrix} \theta_1^{\wti\scrH}(z,x) \\
\phi_1^{\wti\scrH}(z,x) \end{pmatrix}^\top \wti
\Gamma(z)
\begin{pmatrix} \theta_1^{\wti\scrH}(z,x') \\
\phi_1^{\wti\scrH}(z,x') \end{pmatrix} \no \\
&= \sum_{j,k=0}^1 \psi_j(x) \wti \Gamma_{j,k}(z) \psi_k(z), \quad
x<x', \lb{4.47}
\end{align}
where
\begin{equation}
\psi_j(z,x)=\begin{cases} \theta_1^{\wti\scrH}(z,x), & j=0, \\
\phi_1^{\wti\scrH}(z,x), & j=1, \end{cases} \quad (z,x)\in\bbC\times\bbR.
\end{equation}
Next we choose $f_\ell, g_\ell \in C_0^\infty(\bbR)$, $\ell=0,1$, such
that
\begin{equation}
\supp(f_\ell)\subseteq [a,b], \quad \supp(g_\ell)\subseteq [c,d], \quad
a<b<c<d, \; \ell=0,1,
\end{equation}
and introduce the $2\times 2$ matrices
\begin{align}
A(z)&=\left(A_{\ell,m}(z)=
\left\langle \begin{pmatrix} f_\ell \\ 0 \end{pmatrix},
(\wti H-z)^{-1} \begin{pmatrix} g_m \\ 0\end{pmatrix}
\right\rangle_{L^2(\bbR)^2}\right)_{\ell,m=0,1},  \lb{4.50A} \\
& \hspace*{7.5cm}    z\in\rho(\wti H),  \no \\
B(z)&=\left(B_{\ell,j}(z)=\langle f_\ell,\psi_j(z)
\rangle_{L^2(\bbR)}
\right)_{\ell,j=0,1}, \quad z\in\bbC, \lb{4.51} \\
C(z)&=\left(C_{k,m}(z)=\langle \psi_k(z), g_m \rangle_{L^2(\bbR)}
\right)_{k,m=0,1}, \quad z\in\bbC,  \lb{4.52}
\end{align}
where, in obvious notation, $\langle \cdot,\cdot \rangle_{L^2(\bbR)}$
denotes the scalar product in $L^2(\bbR)$ (linear in the second place). In
addition, we let
$z_0\in\rho(\wti H)$ and suppose that $f_\ell$ and $g_\ell$, $\ell=0,1$,
are chosen such that
\begin{equation}
\det (B(z_0))\neq 0, \quad \det(C(z_0))\neq 0.
\end{equation}
Since for fixed $x\in\bbR$, $\psi_j(z,x)$, $j=0,1$, are entire with respect
to $z$, and for fixed $z\in\bbC$, locally absolutely continuous in
$x\in\bbR$, one infers by continuity with respect to $z$ that
\begin{equation}
\det (B(z))\neq 0, \quad \det(C(z))\neq 0, \quad z\in \calU(z_0),
\end{equation}
where $\calU(z_0)\subset \rho(\wti H)$ is a sufficiently small
open neighborhood of $z_0$. Thus, combining \eqref{4.47} and
\eqref{4.50A}--\eqref{4.52} one computes
\begin{equation}
A(z)=B(z)\wti \Gamma(z) C(z), \quad z\in \calU(z_0).
\end{equation}
Since $A$ is analytic on $\calU(z_0)$, $B$ and $C$ are entire and
invertible on $\calU(z_0)$, one concludes that $\wti \Gamma$ is analytic on
$\calU(z_0)$. Since $z_0\in\rho(\wti H)$ was arbitrary, this proves
analyticity of $\wti\Gamma$ on $\rho(\wti H)$.

Since for fixed $x\in\bbR$, $\HwtiF(\cdot,x,\alpha)$ and
$\DwtiF(\cdot,x,\alpha)$ are entire, the claims for $\GopHwti(\cdot,x,x')$
and $\GopDwti(\cdot,x,x')$ are immediate from \eqref{4.34}, \eqref{4.35},
and \eqref{4.45}.
\end{proof}

\subsection{General spectral properties.}\label{s4.3}

In this subsection we recall some of the spectral properties of
$\opDwti$ recorded  in \cite{CGHL04}.

In the following, $\sigma_{\rm
a}(\opDwti), \sigma_{\rm p}(\opDwti), \sigma_{\rm c}(\opDwti),
\sigma_{\rm e}(\opDwti)$, and $\sigma_{\rm r}(\opDwti)$, denote the
approximate point, point, continuous, essential, and residual
spectra of $\wti D$, respectively, while $\pi(\opDwti)$ denotes the
regularity domain and $\rho(\opDwti)$ the resolvent set for
$\opDwti$. Moreover, for $\omega\subset\bbC$, the
complex conjugate of $\omega$ is denoted by
\begin{equation}
\omega^*=\{\ol{\lambda}\in\bbC \,|\, \lambda\in\omega\}. \lb{4.42}
\end{equation}

We begin by noting a result which holds for general $J$-self-adjoint
operators and hence in particular for $\opDwti$ and its unitarily
equivalent $\opHwti$. (Of course, it also applies to the
self-adjoint operators $\hatt D$ and $\hatt H$, cf.\ \eqref{2.14a}.)
\begin{theorem} [\cite{CGHL04}] \lb{t4.9}
Let $\wti D$ be maximally defined as in \eqref{DO}. Then,
\begin{align} \label{SpecPropA}
\sigma(\opDwti)&=\sigma_{\rm p}(\opDwti)\cup\sigma_{\rm c}(\opDwti) \\
&=\sigma_{\rm p}(\opDwti)\cup\sigma_{\rm e}(\opDwti),\\
\sigma_{\rm r}(\opDwti)&=\emptyset,  \\
\sigma_{\rm p}(\opDwti)&=\sigma_{\rm p}(\opDwti^*)^*,  \lb{4.39} \\
\sigma_{\rm a}(\opDwti)&=\sigma(\opDwti), \\
\pi(\opDwti)&=\rho(\opDwti).
\end{align}
\end{theorem}
Remark \ref{r4.1.2} is a crucial ingredient in the proof of the next
result which details spectral properties specific to $\opDwti$ but
which extend to $\wti H$ by the unitary equivalence found in
\eqref{4.6a}.
\begin{theorem} [\cite{CGHL04}] \lb{t4.10}
Let $\wti D$ be maximally defined as in \eqref{DO}. Then,
\begin{align} \label{SpecPropB}
&\sigma(\opDwti)^*=\sigma(\opDwti),\\
&\sigma_{\rm c}(\opDwti)\supseteq \bbR,  \\
&\sigma_{\rm e}(\opDwti)\supseteq \bbR,  \\
&\sigma_{\rm p}(\opDwti)\cap\bbR=\emptyset.
\end{align}
\end{theorem}

Thus, the spectrum for $\opDwti$ is symmetric with respect to
$\bbR$, the continuous spectrum contains $\bbR$, and the point
spectrum is disjoint from $\bbR$. Non-real continuous and essential
spectrum can occur as is seen in the example to follow in which
\emph{crossing spectral arcs} are an essential feature. Contrasted
with this is  the fact that the spectrum for the self-adjoint
operator $\opDhat$ is of course confined to $\bbR$.

To underscore the relevance of the Weyl--Titchmarsh coefficients
$\mwti_{\pm}(z,\gamma)$ for spectral theoretic questions concerning
the non-self-adjoint operators $\wti D$ and $\wti H$, we now present
a calculation analogous to that of Theorem \ref{tWSHhat}. Much more
remains to be done in this context and the remainder of this
subsection offers just a preliminary glimpse at the difficulties
imposed by non-self-adjoint Dirac-type operators.

By analogy with the self-adjoint case discussed in
Subsection~\ref{GMSA}, we denote by $\wti
M(\cdot,\gamma)\in\bbC^{2\times 2}$, $z\in\rho(\wti H)$, the
whole-line Weyl--Titchmarsh $M$-function of the operator $\opHwti$
defined in \eqref{HO} in association with the special case given by
\eqref{HENSA},
\begin{align}
\Mwti(z,\gamma)&= \big( \wti M_{\ell,\ell'}(z,\gamma)\big)_{\ell,\ell'=0,1}
=\frac12[\wti\Gamma(z,\gamma)+
\wti\Gamma(z,\gamma)^\top] = \wti\Gamma(z,\gamma)+
\f{1}{2} \begin{pmatrix} 0 & -1 \\ 1 & 0 \end{pmatrix} \no \\
&=\begin{pmatrix} \frac1{\mwti_-(z,\gamma) -
\mwti_+(z,\gamma)}& \frac12\frac{\mwti_-(z,\gamma) +
\mwti_+(z,\gamma)}{\mwti_-(z,\gamma) -
\mwti_+(z,\gamma)}\\[6pt]
\frac12\frac{\mwti_-(z,\gamma) +
\mwti_+(z,\gamma)}{\mwti_-(z,\gamma) - \mwti_+(z,\gamma)}&
\frac{\mwti_-(z,\gamma)\mwti_+(z,\gamma)}{\mwti_-(z,\gamma) -
\mwti_+(z,\gamma)}
\end{pmatrix}, \quad z\in\rho(\wti H)=\rho(\wti D).  \label{Mwti}
\end{align}
Here, by Remark \ref{r4.4},
$\mwti_\pm(z,\gamma)=\Hwtim_\pm(z,\alpha)=\Dwtim_\pm(z,\beta)$ for
$\beta=\alpha U$, and $\alpha=(\cos(\theta),\sin(\theta))$, $\theta\in
[0,2\pi)$. By \eqref{4.45},
\begin{equation}
\text{$\wti M(\cdot,\gamma)$ is analytic on $\rho(\wti H)=\rho(\wti D)$.}
\end{equation}

Given $\wti M (\cdot,\gamma)$, we introduce the set function
$\wti\Omega(\cdot,\gamma)$ on intervals $(\lambda_1,\lambda_2]\subset
\bbR$, $\lambda_1<\lambda_2$, by
\begin{equation}
\wti\Omega((\lambda_1,\lambda_2],\gamma)=\f{1}{2\pi i}
\lim_{\delta\downarrow 0}\lim_{\varepsilon\downarrow 0}
\int_{\lambda_1+\delta}^{\lambda_2+\delta} d\lambda \,
\big[\wti M(\lambda+i\varepsilon,\gamma)-
\wti M(\lambda-i\varepsilon,\gamma)\big].  \lb{4.48}
\end{equation}
To proceed as in the self-adjoint case in Subsection \ref{GMSA}, we now
make the following set of assumptions.
\begin{hypothesis} \lb{h4.14} Let $[\lambda_1,\lambda_2]\subset\bbR$,
$\lambda_1<\lambda_2$. \\
$(i)$ Suppose no spectral component of $\wti H$ other than
$[\lambda_1,\lambda_2]$ intersects the set
\begin{equation}
\{z\in\bbC\,|\, \lambda_1\leq \Re(z)\leq\lambda_2, \,
0\leq |\Im(z)|\leq \varepsilon_0 \},  \lb{4.50a}
\end{equation}
for some fixed $\varepsilon_0>0$. \\
$(ii)$ Assume that \eqref{4.48} defines a measure on the Borel subsets of
$[\lambda_1,\lambda_2]$. \\
$(iii)$ Suppose that
\begin{align}
\begin{split}
&\varepsilon|\wti M_{\ell,\ell'}(\lambda+i\varepsilon,\ga)|\leq
C(\lambda_1,\lambda_2,\varepsilon_0), \quad \lambda\in
[\lambda_1,\lambda_2], \; 0<\varepsilon\leq\varepsilon_0, \;
\ell,\ell'=0,1,  \\
&\varepsilon
|\wti M_{\ell,\ell'}(\lambda+i\varepsilon,\ga)+
\wti M_{\ell,\ell'}(\lambda-i\varepsilon,\ga)|
\underset{\varepsilon\downarrow 0}{=}\oh(1), \quad \lambda\in
[\lambda_1,\lambda_2], \; \ell,\ell'=0,1.    \lb{4.71A}
\end{split}
\end{align}
\end{hypothesis}

We also use the abbreviation
\begin{equation}
\big( \HwtiT_0(\alpha) f\big)(\lambda) = \int_\bbR dx \,
\HwtiF(\lambda,x,\alpha)^\top f(x),
\quad \lambda \in [\lambda_1,\lambda_2], \; f\in
C^\infty_0(\bbR)^2.
\end{equation}

Analogous definitions and hypotheses apply, of course, to other
parts of the spectrum, assuming one can separate a (complex)
neighborhood of the spectral arc in question from the rest of the
spectrum of $\wti H$ similarly to \eqref{4.50a}. (We note that this
excludes the possibility of crossing spectral arcs, cf.\
\cite{GT06} and Lemma \ref{l4.16}). The extent to which Hypothesis
\ref{h4.14} applies to general $J$-self-adjoint operators studied in
this section is beyond the scope of this paper and will be taken up
elsewhere. Typical examples we have in mind are periodic and certain
classes of quasi-periodic operators $\wti H$, where the spectrum is
known to consist of piecewise analytic arcs.

Given an interval $(\lambda_1,\lambda_2]$ with properties as in
Hypothesis \ref{h4.14}, we define the analog of the spectral
projection \eqref{3.46} in the self-adjoint case, now denoted by
$E_{\wti H}((\lambda_1,\lambda_2])$, associated  with $\wti H$ and
$(\lambda_1,\lambda_2]$, by
\begin{align}
& \langle f, E_{\opHwti}((\lambda_1,\lambda_2]) g\rangle_{L^2(\bbR)^2}
\no \\
& \quad =\lim_{\delta \downarrow 0}\lim_{\varepsilon \downarrow
0}\frac1{2\pi i}\int_{\lambda_1+\delta}^{\lambda_2+\delta}d\lambda \,
\big\langle f,\big[(\opHwti-(\lambda+i\varepsilon))^{-1} \no \\
& \hspace*{4.1cm}
-(\opHwti-(\lambda-i\varepsilon))^{-1}\big] g\big\rangle_{L^2(\bbR)^2},
\quad f,\, g\in C^\infty_0(\bbR)^2.  \lb{4.73}
\end{align}
In the present non-self-adjoint context, it is far from obvious that
$E_{\wti H}((\lambda_1,\lambda_2])$ extends to a bounded operator, let
alone, a bounded projection, on $L^2(\bbR)^2$. A careful study of this
question is again beyond the scope of this paper and hence we introduce the
following hypothesis for now and postpone a detailed discussion
of the properties of $E_{\wti H}((\lambda_1,\lambda_2])$ to a future
investigation:

\begin{hypothesis}\label{h4.15}
Given an interval $(\lambda_1,\lambda_2]$ with properties as in
Hypothesis \ref{h4.14}, we suppose that
$E_{\wti H}((\lambda_1,\lambda_2])$, as defined in \eqref{4.73}, extends
to a bounded projection operator on $L^2(\bbR)^2$.
\end{hypothesis}

We note that Hypotheses \ref{h4.14} and \ref{h4.15} can be verified
in some special cases. For instance, in the case of periodic
Schr\"odinger operators, one can successfully apply Floquet theory
and verify Hypothesis \ref{h4.14} in connection with parts of
spectral arcs which are not intersected by other spectral arcs (cf.\
\cite{GT06}). On the other hand, Hypothesis \ref{h4.15} is known to
fail in the presence of crossings of spectral arcs as shown in
\cite{GT06}. This is also underscored in Lemma \ref{l4.16} in
connection with the simple constant coefficient Dirac-type operator
$\wti D_{q_0}$, which exhibits the crossing of spectral arcs at the
origin. We will return to this circle of ideas elsewhere.

\begin{theorem}\label{t4.11}
Assume Hypotheses \ref{h4.14} and \ref{h4.15} and
$f,\, g\in C_0^\infty(\bbR)^2$. Then,
\begin{equation}\label{4.50}
\langle f, E_{\opHwti}((\lambda_1,\lambda_2])g\rangle_{L^2(\bbR)^2} =
\int_{(\lambda_1,\lambda_2]}\left(\big( \HwtiT_0(\alpha_0)
\ol f\big)(\lambda)\right)^\top d\wti\Omega(\lambda,\gamma_0)
\big(\HwtiT_0(\alpha_0) g\big)(\lambda).
\end{equation}
\end{theorem}
\begin{proof}
For simplicity we will suppress the $\alpha$ (resp., $\gamma$) dependence
of all quantities involved. We closely follow the strategy of proof
employed in connection with Theorem \ref{tWSHhat}.  Then,
\begin{align}
&\langle f, E_{\opHwti}((\lambda_1,\lambda_2])g\rangle_{L^2(\bbR)^2}
\lb{4.75A}  \\
& \quad =\lim_{\delta \downarrow 0}\lim_{\varepsilon \downarrow 0}
\frac1{2\pi i}\int_{\lambda_1+\delta}^{\lambda_2+\delta}d\lambda \,
\big\langle f,\
\big[(\opHwti-(\lambda+i\varepsilon))^{-1}-
(\opHwti-(\lambda-i\varepsilon))^{-1}\big]g \big\rangle_{L^2(\bbR)^2}
\notag \\
& \quad =\lim_{\delta \downarrow 0}\lim_{\varepsilon
\downarrow 0}\frac1{2\pi
i}\int_{\lambda_1+\delta}^{\lambda_2+\delta}d\lambda
\int_{\bbR}dx\int_{\bbR}dx'\big\{
f(x)^*\GopHwti(\lambda+i\varepsilon, x,x')g(x')  \no \\
& \hspace*{6cm} - f(x)^*\GopHwti(\lambda-i\varepsilon,x,x')g(x')\big\} \no
\\ & \quad =\lim_{\delta \downarrow 0}\lim_{\varepsilon \downarrow
0}\frac1{2\pi i}\int_{\lambda_1+\delta}^{\lambda_2+\delta}d\lambda
\int_{\bbR}dx\bigg\{\int_{-\infty}^xdx'
\big\{f(x)^*\GopHwti(\lambda+i\varepsilon,x,x')g(x') \no \\
& \hspace*{6.6cm}  - f(x)^*\GopHwti(\lambda-i\varepsilon,x,x')g(x')\big\}
\no \\
& \hspace*{5.1cm}
+\int_x^{\infty}dx'\big\{f(x)^*\GopHwti(\lambda+i\varepsilon,x,x')g(x')
\no \\ & \hspace*{6.2cm}
-f(x)^*\GopHwti(\lambda-i\varepsilon,x,x')g(x')\big\}\bigg\} \no \\
& \quad
=\lim_{\delta \downarrow 0}\lim_{\varepsilon \downarrow 0}\frac1
{2\pi i}\int_{\lambda_1+\delta}^{\lambda_2+\delta}d\lambda\int_{\bbR}dx
\no \\
&\qquad \times \bigg\{\int_{-\infty}^xdx'
f(x)^*\HwtiF(\lambda,x)\Big[\Gammawti(\lambda+i\varepsilon)
- \Gammawti(\lambda-i\varepsilon) \,
\Big]\HwtiF(\lambda,x')^\top g(x') \no \\
& \hspace*{1.4cm} +\int_x^{\infty}dx'f(x)^*\HwtiF(\lambda,x)
\Big[\Gammawti(\lambda+i\varepsilon)^\top
- \Gammawti(\lambda-i\varepsilon)^\top
\, \Big]\HwtiF(\lambda,x')^\top g(x')\bigg\}.  \no
\end{align}
Here we used conditions \eqref{4.71A} to pass to the last line in
\eqref{4.75A}. By means of the fundamental identity given in \eqref{4.16}, a
calculation shows that
\begin{equation}
\Gammawti(\lambda+i\varepsilon) - \Gammawti(\lambda-i\varepsilon)=
[\Gammawti(\lambda+i\varepsilon) -
\Gammawti(\lambda-i\varepsilon)]^\top,
\end{equation}
and by \eqref{Mwti} that
\begin{equation}
\Mwti(\lambda+i\varepsilon)-\Mwti(\lambda-i\varepsilon)=
\Gammawti(\lambda+i\varepsilon) - \Gammawti(\lambda-i\varepsilon).
\end{equation}
As a consequence,
\begin{align}\label{4.53}
&\langle f, E_{\opHwti}((\lambda_1,\lambda_2])g\rangle_{L^2(\bbR)^2}  \\
&\quad = \int_{(\lambda_1,\lambda_2]} \int_{\bbR}dx
\int_\bbR dx' f(x)^*\HwtiF(\lambda,x) d\wti \Omega (\lambda)
\HwtiF(\lambda,x')^\top g(x'),\notag
\end{align}
from which  \eqref{4.50} then follows.
\end{proof}

As a corollary, we obtain the corresponding result for the operator
$\opDwti$.

\begin{corollary}  \lb{c4.17}
Assume Hypotheses \ref{h4.14} and \ref{h4.15} and $f,\, g\in
C_0^\infty(\bbR)^2$. Then,
\begin{equation}  \lb{4.79}
\langle f, E_{\opDwti}((\lambda_1,\lambda_2])g\rangle_{L^2(\bbR)^2}  =
\int_{(\lambda_1,\lambda_2]}\left(\big( \HwtiT_0(\alpha_0)
\ol{Uf} \big)(\lambda)\right)^\top
d\wti\Omega(\lambda,\gamma_0)\big(\HwtiT_0(\alpha_0) Ug\big)(\lambda).
\end{equation}
\end{corollary}
\begin{proof}
This follows immediately from the observation that
\begin{align}
&\langle f, E_{\opDwti}((\lambda_1,\lambda_2])g\rangle_{L^2(\bbR)^2}  \\
& \quad
=\lim_{\delta \downarrow 0}\lim_{\varepsilon \downarrow 0}\frac1{2\pi
i}\int_{\lambda_1+\delta}^{\lambda_2+\delta}d\lambda \, \big\langle f,
U^{-1}\big[(\opHwti-(\lambda+i\varepsilon))^{-1}  \no \\
&\hspace*{4.6cm} -
(\opHwti-(\lambda-i\varepsilon))^{-1}\big]Ug \big\rangle_{L^2(\bbR)^2}.
\notag
\end{align}
\end{proof}

\subsection{An Example.}  \lb{s4.4}

Next we turn to the calculation of quantities discussed in the
previous section for $\Hwti$ and $\Dwti$ in the case where
$q(x)=q_0\in\bbC\backslash\{0\}$ is constant. (The self-adjoint case
$q_0=0$ has already been discussed in Subsection~\ref{s3.4}.) In this case
we denote $\Hwti$ and $\Dwti$ by $\Hwti_{q_0}$ and $\Dwti_{q_0}$, etc.
Without much exaggeration, the example to follow describes probably the
simplest periodic (even constant coefficient) differential operator with
crossing spectral arcs.

In consideration of this case, we define the function $\Rq(z)$ to be
a function that is analytic with positive imaginary part  on the
split plane
\begin{equation}
\wti\bbP_{q_0}=\bbC\backslash (\bbR\cup\{z\in\bbC\,|\,
z=it,\, t\in[-|q_0|,|q_0|]\})
\end{equation}
such that
\begin{equation}
\Rq(z) = \sqrt{z^2 + |q_0|^2}, \quad z\in\wti\bbP_{q_0}.
\end{equation}
Given this definition, the following conventions are used:
\begin{equation}
\Rq(\ol z)=-\ol{\Rq(z)},\quad z\in\wti\bbP_{q_0},
\end{equation}
\begin{equation}
\Rq(\lambda \pm i0)=\lim_{\varepsilon\downarrow 0}\Rq(\lambda\pm
i\varepsilon)=
\begin{cases}
\pm \sqrt{\lambda^2+|q_0|^2}, & \lambda >0,\\
\mp \sqrt{\lambda^2+|q_0|^2}, & \lambda <0,
\end{cases}\\
\end{equation}
\begin{equation}
\Rq(it \pm 0)=\lim_{\lambda\downarrow 0}\Rq(it \pm \lambda)=
\begin{cases}
\pm\sqrt{|q_0|^2-t^2},  & 0<t\le |q_0|,\\
\mp\sqrt{|q_0|^2-t^2},  & -|q_0|\le t <0.
\end{cases}
\end{equation}

A direct calculation shows for general
$\alpha=(\cos(\theta),\sin(\theta))$, $\theta\in [0,2\pi)$, that
\begin{equation}
\HwtiqPsi_\pm(z,x,\alpha) = a_\pm
\begin{pmatrix}1\\
\frac{\Im(q_0) \mp  \Rq(z)}{iz+\Re(q_0)}
\end{pmatrix}e^{\pm i \Rq(z) x}, \quad z\in\wti\bbP_{q_0}.
\end{equation}
For the the corresponding general $\beta=\alpha
U=[(-1+i)/2] (e^{-i\theta},e^{i\theta})$, $\theta\in [0,2\pi)$,
we  see by direct calculation that
\begin{equation}
\Psi^{\Dhat_{q_0}}_\pm(z,x,\beta)=
b_\pm\begin{pmatrix}1\\ \frac{i}{q_0}(z\pm \wti S_{q_0}(z))\end{pmatrix}
e^{\pm i \wti S_{q_0}(z)x}, \quad z\in\wti\bbP_{q_0}
\end{equation}
for some $b_\pm\in\bbC$, and alternatively that
\begin{align}
\DwtiqPsi_\pm(z,x,\beta)&=U^{-1}\HwtiqPsi_\pm(z,x,\alpha)  \no \\
&=\frac{a_\pm(-1+i)}{2(iz+\Re(q_0))}
\begin{pmatrix}
-z+iq_0 \pm \Rq(z)\\
-z+i\ol q_0 \mp \Rq(z)
\end{pmatrix}e^{\pm i\Rq(z)x}, \\
& \hspace{6.2cm} \quad z\in\wti\bbP_{q_0}  \no
\end{align}
for some $a_\pm\in\bbC$.

In particular, for $\alpha=\alpha_0=(1,0)$ and $z\in\wti\bbP_{q_0}$,
\begin{equation}
\Psi^{\Hhat_{q_0}}_\pm(z,0,\alpha_0) =\calF^{\Hwti_{q_0}}(z,0,\alpha_0)
\begin{pmatrix}1\\
m^{\Hwti_{q_0}}_\pm(z,\gamma_0)\end{pmatrix} =a_\pm \begin{pmatrix}1\\
\frac{\Im(q_0) \mp  \Rq(z)}{iz+\Re(q_0)}
\end{pmatrix}.
\end{equation}
Because $\calF^{\Hwti_{q_0}}(z,0,\alpha_0) = I_2$, we conclude that
$a_\pm=1$ and that
\begin{equation}
\mwti_{q_0,\pm}(z,\gamma_0)=m^{\Hwti_{q_0}}_\pm(z,\alpha_0)
=m^{\Dwti_{q_0}}_\pm(z,\beta_0)= \frac{\Im(q_0)
\mp  \Rq(z)}{iz+\Re(q_0)}, \quad z\in\wti\bbP_{q_0}.
\end{equation}

By \eqref{4.7} we see that
\begin{equation}  \lb{4.68}
\begin{aligned}[t]
&G^{\opHwti_{q_0}}(z,x,x')\\
&=\frac{1}{2\Rq(z)}\begin{pmatrix}iz+\Re(q_0)& \Im(q_0)\mp \Rq(z)\\[2pt]
\Im(q_0)\pm \Rq(z)& iz-\Re(q_0)\end{pmatrix}e^{\pm i\Rq(z)(x'-x)},\\
& \hspace{7.3cm} x\lessgtr x',\ z\in\wti\bbP_{q_0}.
\end{aligned}
\end{equation}
Using  either \eqref{4.6} or the fact that
$G^{\opDwti_{q_0}}(z,x,x')=U^{-1}G^{\opHwti_{q_0}} (z,x,x')U$, we see
that
\begin{equation}
\begin{aligned}[t]
&G^{\opDwti_{q_0}} (z,x,x')=\frac{1}{2\Rq(z)}
\begin{pmatrix}
i\big[z\pm \Rq(z)\big] & q_0 \\[2pt]
\ol{q_0} & i\big[z\mp \Rq(z)\big]
\end{pmatrix}e^{\pm i\Rq(z)(x'-x)},\\
& \hspace{8.35cm} x\lessgtr x',\ z\in\wti\bbP_{q_0}.  \lb{4.69}
\end{aligned}
\end{equation}
The spectra of $\wti{H}_{q_0}$ and $\wti{D}_{q_0}$ are purely
continous and given by
\begin{equation}
\sigma(\wti{H}_{q_0})= \sigma (\wti{D}_{q_0}) =\bbR \cup
\{z\in\bbC\,|\, z=it, \, t\in[-|q_0|,|q_0|]\},
\end{equation}
that is, the spectrum consists of the real axis and the interval from
$-|q_0|$ to $|q_0|$ along the imaginary axis. Since
$q_0\in\bbC\backslash\{0\}$, the origin is a crossing point of the two
spectral arcs.

In contrast to the self-adjoint example $q_0\in\bbC\backslash\{0\}$
discussed in Section~\ref{s3.4}, the potential pole for
$\mwti_{q_0,\pm}(z,\gamma_0)$ given by $z=i\Re(q_0)$ now lies in the
continuous spectrum for $\opHwti$.  For $z=i\Re(q_0)$ and $\Im (q_0)
\ne 0$, the $z$-wave functions for $\Hwti_{q_0}$ are given by
\begin{equation}
\Psi^{\Hwti_{q_0}}(x)= \psi_1(0)\begin{pmatrix}e^{ix\Im (q_0) } \\
2\left(\frac{i\Re (q_0)}{\Im (q_0)}\right)\sin (x\,\Im (q_0) )\end{pmatrix}
+\psi_2(0)\begin{pmatrix}0\\ e^{-ix\Im (q_0)}  \end{pmatrix},
\end{equation}
and for $z=i\Re(q_0)$ and $\Im (q_0) = 0$, the $z$-wave functions for
$\Hwti_{q_0}$ are given by
\begin{equation}
\Psi(x)= \psi_1(0)\begin{pmatrix}1\\-2ix \, \Re (q_0) \end{pmatrix} +
\psi_2(0)\begin{pmatrix}0\\1 \end{pmatrix}.
\end{equation}
Consequently, $z= i\Re (q_0)$ is  neither an eigenvalue of
$\opHwti_{q_0}$ nor an eigenvalue of $\opHwti_{q_0,\pm} (\alpha_0)$.
\subsection{Nonspectrality}  \lb{s4.5}

The principal result of this subsection illustrates that for all the
similarities developed thus far, $\opDhat$ and $\opDwti$ bear the following
stark difference: $\opDhat$, being self-adjoint, is always a
spectral operator of scalar type in the sense of Dunford and
Schwartz while $\opDwti$ cannot be expected to be a spectral
operator whenever there are crossing spectral arcs in the spectrum
of $\opDwti$.

In the case of periodic Schr\"odinger operators, this result has recently
been proved in \cite{GT06}. Here we confine ourselves to a study of the
constant coefficient operator $\wti{D}_{q_0}$ but on the basis of
\cite{GT06} it is natural to expect this result extends to all periodic
Dirac-type operators $\opDwti$ with crossing spectral arcs.

Applying Corollary \ref{c4.17} to the concrete example
$q(x)=q_0\in\bbC\backslash\{0\}$ treated in the previous subsection, one
can rewrite \eqref{4.79} to obtain
\begin{align}
& \langle f,
E_{\opDwti_{q_0}}((\lambda_1,\lambda_2))g\rangle_{L^2(\bbR)^2}  \\ &
\quad =\int_{(\lambda_1,\lambda_2)}\left(\big(T_0^{\wti\scrH_{q_0}}(\alpha)
\ol{U f}\big)(\lambda)\right)^\top d\wti\Omega_{q_0}(\lambda,\gamma)
\big(T_0^{\wti\scrH_{q_0}}(\alpha) U g\big)(\lambda)  \no \\
& \quad = \f{1}{2} \int_{\lambda_1}^{\lambda_2} d\lambda \, \bigg[
\begin{pmatrix} \hat f_1(\sqrt{\lambda^2+|q_0|^2})  \\
\hat f_2(\sqrt{\lambda^2+|q_0|^2})\end{pmatrix}^*  \no \\
& \hspace*{1.95cm} \times
\begin{pmatrix} \big[\lambda/\sqrt{\lambda^2+|q_0|^2}\,\big]+1 &
-iq_0/\sqrt{\lambda^2+|q_0|^2} \\ -i\ol{q_0}/\sqrt{\lambda^2+|q_0|^2} &
\big[\lambda/\sqrt{\lambda^2+|q_0|^2}\,\big]-1
\end{pmatrix}
\begin{pmatrix} \hat g_1(\sqrt{\lambda^2+|q_0|^2})
\\ \hat g_2(\sqrt{\lambda^2+|q_0|^2})\end{pmatrix} \no \\
& \hspace*{1.5cm} +\begin{pmatrix} \check f_1(\sqrt{\lambda^2+|q_0|^2})
\\ \check f_2(\sqrt{\lambda^2+|q_0|^2})\end{pmatrix}^*
\begin{pmatrix} \big[\lambda/\sqrt{\lambda^2+|q_0|^2}\,\big]-1 &
-iq_0/\sqrt{\lambda^2+|q_0|^2} \\ -i\ol{q_0}/\sqrt{\lambda^2+|q_0|^2} &
\big[\lambda/\sqrt{\lambda^2+|q_0|^2}\,\big]+1 \end{pmatrix}  \no \\
& \hspace*{2cm} \times
\begin{pmatrix} \check g_1(\sqrt{\lambda^2+|q_0|^2})
\\ \check g_2(\sqrt{\lambda^2+|q_0|^2})\end{pmatrix}\bigg],
\quad 0<\lambda_1<\lambda_2, \; f,g \in C_0^\infty(\bbR)^2,  \no
\end{align}
where
\begin{align}
& \hat h (p)=\f{1}{(2\pi)^{1/2}} \int_{\bbR} dx\, e^{ipx} h(x), \quad
\check h (p)=\f{1}{(2\pi)^{1/2}} \int_{\bbR} dx\, e^{-ipx} h(x), \\
& \hspace*{6.95cm} p\in\bbR, \; h\in C_0^\infty(\bbR). \no
\end{align}

We note that the spectrum of $\opDwti_{q_0}$ is purely continuous and so
the distinction between the intervals $(\lambda_1,\lambda_2]$ and
$(\lambda_1, \lambda_2)$ becomes irrelevant throughout this subsection.

In the following, $\calB(\calH)$ denotes the Banach space of bounded linear
operators in a Hilbert space $\calH$.

\begin{lemma}\label{l4.16}
Let $[\lambda_1,\lambda_2]\subset (0,\infty)$. Then,
\begin{equation}
\lim_{\lambda_1\downarrow 0} \|E_{\opDwti_{q_0}}
((\lambda_1,\lambda_2))\|_{\calB(L^2(\bbR)^2)}
=\infty.   \lb{4.72}
\end{equation}
\end{lemma}
\begin{proof}
We choose $f,g$ of the form $f=(h, 0)^\top$, $g=(0, h)^\top$, $h\in
C_0^\infty(\bbR)$. Then,
\begin{align}
& \langle f, E_{\opDwti_{q_0}}((\lambda_1,\lambda_2))g\rangle_{L^2(\bbR)^2}
\\
& \quad =
-\int_{\lambda_1}^{\lambda_2} d\lambda \f{iq_0}{2\sqrt{\lambda^2+|q_0|^2}}
\big[|\hat h(\sqrt{\lambda^2+|q_0|^2})|^2 + |\check
h(\sqrt{\lambda^2+|q_0|^2})|^2 \,\big].
\end{align}
The change of variables
\begin{equation}
\mu=\sqrt{\lambda^2+|q_0|^2} \geq |q_0|, \quad d\lambda=\f{\mu
d\mu}{\sqrt{\mu^2-|q_0|^2}}
\end{equation}
then yields
\begin{equation}
\langle f, E_{\opDwti_{q_0}}((\lambda_1,\lambda_2))g\rangle_{L^2(\bbR)^2}
=-\f{iq_0}{2}
\int_{\sqrt{\lambda_1^2+|q_0|^2}}^{\sqrt{\lambda_2^2+|q_0|^2}}
\f{d\mu}{\sqrt{\mu^2-|q_0|^2}} \big[|\hat h(\mu)|^2
+|\check h(\mu)|^2  \,\big].  \lb{4.76}
\end{equation}
It suffices to study the first term on the right-hand side of
\eqref{4.76}. (The second term is handled in exactly the same manner.)
For this purpose we now introduce in $L^2(\bbR; d\mu)$ the maximally
defined operator $T(\lambda_1,\lambda_2)$ of multiplication by the function
\begin{equation}
t((\lambda_1,\lambda_2), \mu)=\f{1}{\sqrt{\mu^2-|q_0|^2}}
\chi_{[\sqrt{\lambda_1^2+|q_0|^2},\sqrt{\lambda_2^2+|q_0|^2}]}(\mu), \quad
\mu\in\bbR,
\end{equation}
where $\chi_\omega$ denotes the characteristic function of the set
$\omega\subset\bbR$. We recall that
\begin{equation}
\|T(\lambda_1,\lambda_2)\|_{\calB(L^2(\bbR;d\mu))}
=\|t((\lambda_1,\lambda_2),\cdot)\|_{L^\infty(\bbR;d\mu)}
\end{equation}
(cf.\ \cite[p.\ 51--54]{We80}).

Next, we note that
\begin{equation}
\int_{\sqrt{\lambda_1^2+|q_0|^2}}^{\sqrt{\lambda_2^2+|q_0|^2}}
\f{d\mu}{\sqrt{\mu^2-|q_0|^2}} |\hat h(\mu)|^2 =
\|T(\lambda_1,\lambda_2) \hat h\|_{L^2(\bbR;d\mu)}^2.
\end{equation}
Thus, as long as $0<\lambda_1<\lambda_2$, one infers that $\mu>|q_0|$ and
hence that
\begin{equation}
\|T(\lambda_1,\lambda_2)\|_{\calB(L^2(\bbR;d\mu))}
=\|t((\lambda_1,\lambda_2),\cdot)\|_{L^\infty(\bbR;d\mu)}<\infty,
\quad 0<\lambda_1<\lambda_2.
\end{equation}
However, since $\mu\downarrow q_0$ as $\lambda\downarrow 0$, one obtains
\begin{equation}
\|T(\lambda_1,\lambda_2)\|_{\calB(L^2(\bbR;d\mu))}
=\|t((\lambda_1,\lambda_2),\cdot)\|_{L^\infty(\bbR;d\mu)}\uparrow\infty
\, \text{ as $\lambda_1\downarrow 0$}.  \lb{4.107}
\end{equation}
Since $\|\hat h\|_{L^2(\bbR)}=\|\check h\|_{L^2(\bbR)}=\|h\|_{L^2(\bbR)}$,
there exists a $C>0$ such that
\begin{equation}
\int_{\sqrt{\lambda_1^2+|q_0|^2}}^{\sqrt{\lambda_2^2+|q_0|^2}}
\f{d\mu}{\sqrt{\mu^2-|q_0|^2}} |\hat h(\mu)|^2 \leq C \|h\|_{L^2(\bbR)}^2
\end{equation}
for all $h\in C^\infty_0(\bbR)$ if and only if
$T(\lambda_1,\lambda_2)\in\calB(L^2(\bbR; d\mu))$ and hence if and only if
$t((\lambda_1,\lambda_2),\cdot)\in L^\infty(\bbR; d\mu)$. The blowup in
\eqref{4.107} then shows that \eqref{4.72} holds.
\end{proof}

Thus, the crossing of spectral arcs at the point $\lambda=0$
prevents the operator $\opDwti_{q_0}$ (and hence $\opHwti_{q_0}$) to
have a uniformly bounded family of spectral projections. This is
remarkable since the corresponding Green's matrices \eqref{4.68} and
\eqref{4.69} exhibit no singularity at $z=0$.

\bigskip
\noindent {\bf Acknowledgements.}
We would like to thank the organizers for creating a
great scientific atmosphere and for their extraordinary efforts and
hospitality during the 2005 UAB International Conference on Differential
Equations and Mathematical Physics. 

We also gratefully acknowledge pertinent hints to the literature by E.\
Brian Davies and Alexander Sakhnovich and a critical reading of our
manuscript by Kwang Shin.


\end{document}